\documentclass{aptpub}
\authornames{K. AGARWAL AND V. KAVITHA} 
\shorttitle{Branching processes using Stochastic Approximation} 

\usepackage{graphicx}
\usepackage{enumitem}
\usepackage{xcolor}
\usepackage{nicefrac}

\usepackage{bm}
\usepackage{wrapfig}
\usepackage{amsmath}
\usepackage{amssymb}

\usepackage{imsart}
\usepackage{hyperref}

\begin{document}

\title{New Results in branching processes using stochastic approximation
} 

\authorone[IEOR, IIT Bombay]{KHUSHBOO AGARWAL, VEERARUNA KAVITHA} 

\emailone{agarwal.khushboo@iitb.ac.in and vkavitha@iitb.ac.in} 

\begin{abstract}
We consider a broad class of continuous-time two-type population size-dependent Markov Branching Processes. The offspring distribution can depend on the current (alive) and total (dead and alive) populations. Using stochastic approximation techniques, we show that the time-asymptotic proportion of the populations  either converges to the equilibrium points or infinitely often enters every neighbourhood and exits some neighbourhood of a saddle point of an appropriate ordinary differential equation with a certain probability (almost surely for the process with attack and proportion-dependent branching process). The result holds under finite second-moment conditions. We also show that certain normalized trajectories of the embedded chain almost surely converge to the solution of the ordinary differential equation uniformly over any finite time window  as time progresses. In addition to extending the analysis of several existing BPs, we analyze two new variants: BP with attack and acquisition, and BP with proportion-dependent offsprings. Using these results, we study competition in viral markets and fake news control on online social networks. 
\end{abstract}

\maketitle





\newcommand{\ifnonauto}[2]{#1}

\newcommand{\DetailK}[1]{} 

\newcommand{\tcprocess}{total-current population-dependent BP }
\newcommand{\tcprocessnospace}{total-current population-dependent BP}
\newcommand{\hide}[1]{}

\newcommand{\psiL}{\psi^c_{\mbox{\tiny o}}}
\newcommand{\psiaL}{\psi^a_{\mbox{\tiny o}}}

\newcommand{\bL}{\beta^c_{\mbox{\tiny o}}}
\newcommand{\baL}{\beta^a_{\mbox{\tiny o}}}

\newcommand{\tL}{\theta^c_{\mbox{\tiny o}}}
\newcommand{\taL}{\theta^a_{\mbox{\tiny o}}}

\newcommand{\N}{\mathcal{N}}
\newcommand{\I}{I_{\theta/\psi}}
\newcommand{\eop}{\hfill{$\square$}}
\newcommand{\nto}{\nrightarrow}
\newcommand{\Om}{\Phi}
\newcommand{\om}{\phi}
 \newcommand{\ups}{{\mbox{\small ${\Upsilon}$}}}
 \newcommand{\Ups}{{\mathbf \Upsilon}}
\newcommand{\Bin}{{\cal B}}
\newcommand{\tp }{\tau^+}
\newcommand{\tm}{\tau^-}
\newcommand{\up}{\uparrow}
\newcommand{\offs}{\Gamma}
\newcommand{\down}{\hspace{.15mm}\downarrow}
\newcommand{\beq}{\begin{eqnarray*}}
\newcommand{\eeq}{\end{eqnarray*}}
\newcommand{\Cx}{C^x}
\newcommand{\Cy}{C^y}
\newcommand{\cx}{c^x}
\newcommand{\cy}{c^y}
\newcommand{\cM}{{\cal E}}
\newcommand{\bcx}{{\overline c^x}}
\newcommand{\bcy}{{\overline c^y}}
\newcommand{\bax}{{\overline a^x}}
\newcommand{\Ax}{A^x}
\newcommand{\Ay}{A^y}
\newcommand{\ax}{a^x}
\newcommand{\ay}{a^y}
\newcommand{\wm}{\widetilde{e}}
\newcommand{\cA}{{\cal A}}
\newcommand{\cR}{{\cal S}}
\newcommand{\cRo}{{\cal R}_e}
\newcommand{\tc}{\theta^c}
\newcommand{\ta}{\theta^a}
\newcommand{\pc}{\psi^c}
\newcommand{\pa}{\psi^a}
\newcommand{\Tc}{\Theta^c}
\newcommand{\Ta}{\Theta^a}
\newcommand{\Pc}{\Psi^c}
\newcommand{\Pa}{\Psi^a}
\newcommand{\bc}{\beta^c}
\newcommand{\ba}{\beta^a}
\newcommand{\Bc}{\mathrm{B}^c}
\newcommand{\Ba}{\mathrm{B}^a}
\newcommand{\Beta}{\mathrm{B}}
\newcommand{\dist}{ d_{st}}
\newcommand{\st}{{\cal T}}
\newcommand{\bstar}{\beta^*}
\newcommand{\minf}{m^\infty}
\newcommand{\ueps}{\overline{\varepsilon} }
\newcommand{\leps}{\underline{\varepsilon}}
\newcommand{\polya}{\mbox{P\'{o}lya} }
\newcommand{\q}{\mathbf{q}}

\newcommand{\ga}{\mathbf{g}}
\newcommand{\gna}{\bm{\varrho}}

\newcommand{\cS}{\mathcal{D}_b}
\newcommand{\cD}{\mathcal{D}}
\newcommand{\cB}{\mathcal{B}}

\newcommand{\bpam}{e}
\newcommand{\bpaum}{\overline{e}}
\newcommand{\bpalm}{\underline{e}}
\newcommand{\propum}{\overline{m}}
\newcommand{\proplm}{\underline{m}}
\newcommand{\propwm}{m^\infty}
\newcommand{\propm}{m}
\renewcommand\thefootnote{\arabic{footnote}}

\newcommand\overlinebelow[1]{\stackunder[1.2pt]{$#1$}{\rule{1.2ex}{.075ex}}}
\def\lc{\left\lceil}   
\def\rc{\right\rceil}

\section{Introduction}
One considers the study of growth patterns and limit proportions to analyze Markov chains that are predominantly transient, like branching processes (BPs) under the super-critical regime (for example, \cite{athreya2004branching, kesten1967limit}). This paper precisely investigates the time-asymptotic proportion of population types for a general class of continuous-time two-type population size-dependent Markov BPs. The offspring depends on the current (alive) as well as the total (alive and dead) populations, and can also be negative to model attack (removal of offspring of another type). We analyze such \textit{total-current population-dependent BPs}  in what we call  \textit{throughout super-critical regime} - the expected number of offspring produced by any individual is strictly greater than one, for all population sizes. We will refer to the proportion of the current population size (of one of the types) as the proportion and the time-asymptotic proportion as the limit proportion.

The literature mainly considers offspring that depend only on the current population; such models are essential in several biological applications (for example, \cite{klebaner1984population, yakovlev2009relative}). Recently, authors in  \cite{agarwal2022saturated, hautphenne2022fluid} introduced total-population dependent BPs; however, both papers analyze the BPs which shift from the super-to-sub critical regime, while we are interested in throughout-super-critical BPs. To the best of our knowledge, no other work considers such total-population dependency.

The importance of limit proportions is discussed in various papers, for example, \cite{ranbir2019decomposable, jagers1969proportions, klebaner1989geometric,  klein1980multitype} and several others. Further, they are crucial objects for the analysis of many applications. For example, authors in \cite{kapsikar2020controlling} design a warning mechanism robust against fake news propagation, where the control depends on the proportion of posts marked as fake. In \cite{agarwal2021co}, we study the relative visibility of advertisement posts defined in terms of the limit proportion of unread copies of posts shared by competing content providers. The limit proportions in prey-predator BP of \cite{coffey1991galton} denote
the proportions in which preys and predators co-survive (if at all).

To analyze proportions, it is sufficient to study the embedded chain of the underlying BP. This study is derived using stochastic approximation (SA) techniques (e.g., \cite{kushner2003stochastic}); we have previously used such an amalgam of SA-based methods in BPs in \cite{agarwal2021co, agarwal2022saturated, kapsikar2020controlling}. In this paper, we include a notion of hovering around saddle points and prove that the sets of attractors and saddle points of an autonomous, non-smooth ordinary differential equation (ODE) almost surely describe the limit proportion. In fact, we prove that the limit set of a single-dimensional ODE suffices. We also prove that the ODE solution approximates certain normalized trajectories of the current and total population sizes over any finite time window. 

Previously, SA based approach has been used in the \polya urn (stochastic process closely related to BPs) literature to investigate limit proportions of the balls of a specific colour (see, for example, \cite{athreya1968embedding, higueras2006central, janson2004functional, arthur1987non}). However, the urn-based literature majorly deals with non-extinction scenarios and considers dependency on the current number of balls (not total) in the urn. Further, to the best of our knowledge, no finite time approximation trajectories exist for \polya urn-based models. Furthermore, we also introduce and analyze `BP with attack', where deletion of offspring (attack) from a population type and addition of the same to the other type (acquisition) occurs, in addition to the production of offspring of own type. 
Furthermore, we also study the new `proportion-dependent BP' where the population dependent mean offsprings are proportion-dependent in the limit. Thus, our paper significantly generalizes the models not only in the BP literature but also in the \polya urn literature by including (total and current) population dependency and negative offspring. We provide a more extensive comparison to the existing results in Section \ref{sec_survey}; extensions of the existing results are also provided.

\textbf{Organization:} The main result is provided in Section \ref{sec_probdesc_mainresult} and proved in Section \ref{sec_proof_thrm1}. The ODE analysis is derived in Section \ref{procedureAR}, while BP with attack and proportion-dependent BP, along with their respective applications are in Section \ref{sec_BPA} and \ref{sec_prop_BP}. Section \ref{sec_6} discusses numerical examples for finite time approximation.

\textbf{Notations:} For convenience, we refer the random variable and the corresponding sequence by the same symbol when the context is clear, for example, $\Ups_n$. We abbreviate infinitely often as i.o. and almost surely as a.s. We also use acronyms like BP, SA and ODE defined in the introduction. 
For any function $f$ and time $\tau$, let $f(\tau^-) := \lim_{t \uparrow \tau} f(t)$ and $f(\tau^+) := \lim_{t \downarrow \tau} f(t)$.

\subsection{Problem description}\label{sec_prob}
Consider $x$ and $y$-types of populations, and  
let $\cx_0, \cy_0$ be their respective initial sizes. The lifetime of any individual of any type is exponentially distributed with parameter $0 < \lambda < \infty$ (i.e., \textit{we consider Markovian BPs}). The time instance at which an individual completes its lifetime is referred to as its `death' time.   

Let $\Cx(t), \Cy(t)$ be the \textit{current population} and $\Ax(t), \Ay(t)$ be the \textit{total population} sizes at time $t$. 
Define $\Om(t) := (\Cx(t), \Cy(t), \Ax(t), \Ay(t))$ 
and observe 
$(\Ax(0), \Ay(0)) = (\cx_0, \cy_0)$. Let $\tau$ be the death time of any individual. Let $\offs_{ij}(\Om(\tau^-))$, with  $i, j \in \{x, y\}$, be integer-valued random variables representing $j$-type offspring produced by an $i$-type parent, conditioned on the sigma algebra $\sigma\{\Om(\tau^-)\}$. 
Basically, when $\Om(\tau^-) = \om$, the random offspring are represented by $\offs_{ij}(\om)$ for each $i, j$.
When an individual of $i$-type dies, the sizes of $i$ and $j$-type populations change by $\offs_{ii}(\Om(\tau^-))$ and $\offs_{ij}(\Om(\tau^-))$ respectively\footnote{For each $i, j$, the distribution of $\offs_{ij}(\Om(\tau^-))$ depends on the population size $\Om(\tau^-)$, and not on the value of the epoch, $\tau$.}. Further, the current size (not the total size) of $i$-type reduces by $1$ due to death. The dynamics can then be written as follows, when an $i$-type parent dies, for $i, j \in \{x, y\}$ and $j \neq i$:
\begin{equation}\label{evolve_x_up_time}
\begin{aligned}
C^i(\tau^+) &= C^i(\tau^-)  + \offs_{ii}(\Om(\tau^-)) - 1, \ \ \  A^i(\tau^+) = A^i(\tau^-)  + \offs_{ii}(\Om(\tau^-)), 
\\
C^j(\tau^+) &= C^j(\tau^-) + \offs_{ij}(\Om(\tau^-)), \ \ \ A^j(\tau^+) = A^j(\tau^-) + \offs_{ij}(\Om(\tau^-)).
\end{aligned}
\end{equation}

We consider a significantly generic framework to study \tcprocessnospace, which includes `attack+acquisition' (acquired individuals change their type);  negative (valued) offspring are used to model such attacks. 

In any BP, the expected/mean offspring plays a determining role in the growth of any population. In this paper, we are keen to analyze the super-critical\footnote{See \cite{athreya1968some, athreya2004branching} for an introduction to super-critical population-independent BPs.} variant of total-current population-dependent BPs, which we define formally in the next few lines.
Let $\om = (\cx, \cy, \ax, \ay)$ be a realisation of the random vector $\Om$. Let \textit{$m_{ij} (\om) := E[ \offs_{ij} (\om) ]$ for $i, j \in \{x, y\}$ represent the conditional expectation of the number of offspring, conditioned on $\om$}; we refer these as mean functions and $M(\om) := [m_{ij}(\om)]$ as mean matrix. Then, any BP which satisfies $m_{ix}(\om) + m_{iy}(\om) > 1$ for each $\om$ and $i\in\{x,y\}$ is called \underline{throughout-super-critical BP}. We assume the following for the random number of offspring conditioned on $\om$, which also ensures such super-criticality:
\begin{enumerate}[label=\textbf{A.\arabic*}, ref=\textbf{A.\arabic*}]
    \item \label{a1} There exist two integrable random variables $\overline{\offs}$ and $\underline{\offs}$ which bound the random offspring as: $\underline{\offs} \leq \offs_{ix}(\om) + \offs_{iy}(\om) \leq \overline{\offs}$ a.s., for each $\om$. Also,  $E[\overline{\offs}^2] < \infty$ and $E[\underline{\offs}] > 1$. Further, $\offs_{ii}(\om) \geq 0$ a.s., for each $i, \om$.
\end{enumerate}
Like the population-independent counterparts, the \tcprocess satisfying \ref{a1} also exhibits \textit{dichotomy}: the sum current population, $S^c(t) := \Cx(t) + \Cy(t)$ either explodes (i.e., $S^c(t) \to \infty$ as $t \to \infty$) exponentially at a rate at least $\lambda(E[\underline{\offs}]-1)$ or gets extinct ($S^c(t) = 0$ for all $t \geq t_e$ where $t_e < \infty$) a.s., by Lemma \ref{lemma_sum_pop} in Appendix \ref{appendix_prelim}. Now, our aim is two-fold: (i) to evaluate the limit proportion, $\lim_{t \to \infty} \frac{\Cx(t)}{\Cx(t) + \Cy(t)}$ in non-extinction paths, and (ii) to derive the deterministic approximate trajectories for the underlying BP.

\section{Main result}\label{sec_probdesc_mainresult}
 
When one considers a process which explodes with time, like a typical BP, it is a common practice to scale the process appropriately such that the scaled process converges to a  finite limit; this enables the asymptotic study of the rate of explosion,   proportions of various components of the process, etc. Further, since we are primarily interested in studying limit proportion, it suffices to analyze the embedded process (discrete-time chain defined at death instances). It is important to observe here that such an embedded process is very different from a corresponding BP in discrete-time.

Let $\tau_n$ be the time at which $n$-th individual dies.
Consider any $n \geq 1$. Let $\Om_n := (\Cx_n, \Cy_n, \Ax_n, \Ay_n)$ be the individual (current and total) populations and $S^c_n :=  \Cx_n + \Cy_n$ be the sum current population, both immediately after $\tau_n$, e.g., $\Cx_n = \Cx(\tau_n^+)$. The current population can get extinct, and thus let $\nu_e := \inf \{n : S_n^c = 0\}$ be the extinction epoch,  
with the usual convention that $\nu_e = \infty$, when $S_n^c > 0$ for all $n$. \textit{For the sake of completion, define $\Om_n := \Om_{\nu_e}$ and $\tau_{n} :=\tau_{\nu_e}$, for all $n \geq \nu_e$, when $\nu_e < \infty$.}

Further, define  $\Pc_n := S^c_{n}/n$ and $\Tc_n := \Cx_{n}/n$. Analogously, define  $\Pa_n$ and $\Ta_n$ for the total population, with $S^a_n := \Ax_n + \Ay_n$.  Let $\Ups_n := (\Pc_n, \Tc_n, \Pa_n, \Ta_n)$, and let $\Ups_0 := (s_0^c, \cx_0, s_0^c, \cx_0)$ denote the initial population, where the initial sum population, $s_0^a = s_0^c := \cx_0 + \cy_0$. \textit{Define $\Bc_n := \Tc_n/\Pc_n  = \Cx_n/S_n^c$ as the proportion of $x$-type population among current population}; observe that conditioned on $\Om_n$, the probability of $x$-type individual dying before others is given by $\Bc_n$ in Markovian BPs. Let $\ups := (\pc, \tc, \pa, \ta)$ be a realisation of $\Ups$ and $\bc = \cx/(\cx+\cy) = \tc/\pc$ be that of $\Bc$. 

In the literature, it has been a common practice to assume that the mean matrix converges to a constant matrix for studying (current) population-dependent BPs (\cite{jagers1997coupling, klebaner1984population, klebaner1989geometric}) and they assume convergence at a certain rate (as in \ref{a2} given below). \textit{We extend such work by allowing our total-current population-dependent mean functions, $m_{ij}(\om)$, to converge to proportion-dependent mean functions, $\minf_{ij}(\bc(\ups))$ (which can further be discontinuous)}, while still using the similar convergence criterion. In other words, the limit mean matrix in our case can depend on the proportion.
\begin{enumerate}[label=\textbf{A.\arabic*}, ref=\textbf{A.\arabic*}]
\setcounter{enumi}{1}
    \item Define $\bc(\ups) := \tc/\pc = \cx/s^c$. As sum current population, $s^c \to \infty$:   \label{a2}
    \begin{align*}
    |m_{ij}(\om) - \minf_{ij}(\bc(\ups))| \leq \frac{1}{(s^c)^\alpha}, \mbox{ for each }i, j \in \{x, y\}, \mbox{ for some } \alpha > 0.
    \end{align*}
    Further, assume the existence of finite $\leps$ and $\ueps$ defined below:
\begin{align}\label{eqn_bound_m_general}
    \leps := \inf_{\bc}\{\minf_{ij}(\bc): i, j \in \{x, y\}\} \mbox{ and }\ueps := \sup_{\bc}\{\minf_{ij}(\bc): i, j \in \{x, y\}\}.
\end{align}
\end{enumerate}
Under \ref{a1}-\ref{a2}, we analyze the ratios $\Ups_n$ using SA techniques; specifically, using the solutions of the following ODE:
\begin{align}\label{eqn_ODE}
\dot{\ups} &= \ga(\ups) = \mathbf{h}(\bc)1_{\{\pc > 0\}} - \ups, \mbox{ where } \mathbf{h}(\bc) := (h_{\psi}^{c},  h_{\theta}^{c},  h_{\psi}^{a},  h_{\theta}^{a}), \mbox{ with} \nonumber \\
        h_{\psi}^{c}(\bc) &= \bc \bigg(\minf_{xx}(\bc) +     \minf_{xy}(\bc)\bigg) + (1-\bc)\bigg(\minf_{yy}(\bc) + \minf_{yx}(\bc)\bigg) - 1,  \nonumber\\
        h_{\theta}^{c}(\bc)  &= \bc \bigg(\minf_{xx}(\bc) - 1\bigg) + (1-\bc) \minf_{yx}(\bc), \\
        h_{\psi}^{a}( \bc) &=  \bc \bigg(\minf_{xx}(\bc) + \minf_{xy}(\bc) \bigg) + (1-\bc)\bigg(\minf_{yy}(\bc) + \minf_{yx}(\bc) \bigg)   \mbox{ \normalsize  and}  \nonumber  \\
        h_{\theta}^{a}( \bc)  &= \bc \minf_{xx}(\bc) + (1-\bc) \minf_{yx}(\bc).  \nonumber
\end{align}Given that the above ODE is autonomous and non-smooth (the right hand side is discontinuous), we next assume the existence of the  unique solution in extended sense\ifnonauto{.}{{\color{red}, which is continuous with respect to (w.r.t.) initial conditions.} } 

\begin{definition}\label{defn_solution}
A function $\ups(\cdot)$ is said to be an \underline{extended solution of ODE \eqref{eqn_ODE}} if it is absolutely continuous, and satisfies the equation \eqref{eqn_ODE} for almost all $t \geq 0$.
\end{definition}
\begin{enumerate}[label=\textbf{A.\arabic*}, ref=\textbf{A.\arabic*}]
\setcounter{enumi}{2}
    \item There exists a unique solution $\ups(\cdot)$ for ODE \eqref{eqn_ODE} in the extended sense  over any bounded interval. \ifnonauto{}{{\color{red}Further, assume that the solution is continuous w.r.t. initial conditions, except $0$ such that  for each $T < \infty$, 
    \begin{align}\label{eqn_continuity}
        \sup_{0 \leq t \leq T} d(\ups_{(n)}(t), \ups_{(\infty)}(t)) \to 0, \mbox{ if } \ups_{(n)}(0) \to \ups_{(\infty)}(0), \mbox{ as }n \to \infty,
    \end{align}
    where $\ups_{(n)}(\cdot)$ is the unique solution with initial condition $\ups_{(n)}(0)$, for all $n \leq \infty$.}} \label{a3}
\end{enumerate}
Assumption \ref{a3} is immediately satisfied by standard results in ODEs if $\minf_{ij}(\cdot)$ are Lipschitz continuous and if there was no indicator, $1_{\{\pc > 0\}}$ (see \cite[Theorem 1, sub-section 1.4]{piccini1984ordinary}). We prove the same for ODE \eqref{eqn_ODE} also when  $\minf_{ij}(\cdot)$ are discontinuous and under certain conditions in Theorem \ref{thrm_attractors_beta} in Section \ref{procedureAR}; such discontinuous functions are typical for BPs with attack.
We now recall the definitions of asymptotically stable and  saddle points for autonomous ODE (see \cite{{piccini1984ordinary}}), that facilitates the desired a.s. convergence of $(\Ups_n)$ - 
the definitions are provided to suit our purpose and  can also be applied for the cases with extended solutions of ODE as in Definition \ref{defn_solution}.
\begin{definition}\label{defn_attractors}
(a) Define open ball, $N_\epsilon(\cA) := \{\ups: d(\ups, \cA) < \epsilon\}$.
A set ${\cA}$ is said to be an \underline{attractor or asymptotically stable} and $\cD_\cA$ is said it to be the domain of attraction for  an autonomous ODE
if: (i) the set $\cA$ is stable, i.e., for any $\epsilon > 0$, there exists a $\delta > 0$ such that the solution 
 of the ODE $\ups(t) \in N_\epsilon(\cA)$ for every $t > 0$, if initial condition $\ups(0) \in N_\delta(\cA)$, and (ii) the solution $\ups(t) \to {\cA}$ as $t \to \infty$, if $\ups(0) \in \cD_\cA$.  
 
\noindent (b) Let $\cA^\complement$ be the complement of $\cA$.
A set ${\cR} \subset \cA^\complement$ is said to be \underline{saddle set} if: (i) every $\ups^* \in \cR $ is an equilibrium point, i.e., $\ga(\ups^*) = 0$, and (ii) 
there exists ${\cal D}_S$ such that $\ups(t) \stackrel{t \to \infty}{\longrightarrow} \cA $   for some $\ups(0) \in \cR^\complement \cap {\cal D}_S$ and   $\ups(t) \stackrel{t \to \infty}{\longrightarrow} \cR $  for some other $\ups(0) \in \cR^\complement \cap {\cal D}_S$. 
\end{definition}

By virtue of ODE structure in \eqref{eqn_ODE}, we will see that the saddle points in our case are attracted exponentially to $\cR$ along a particular affine sub-space, and to $\cA$ in the remaining space as given below (we will prove this in Theorem \ref{thrm_attractors_beta} of Section \ref{procedureAR}).

\begin{definition}\label{defn_q_as}
Any non-zero $\ups^*  \in \cR$  is said to be (quasi) \underline{q-attractor} if (i) for any  $\ups(0) \in {\mathbb S}(\ups) := \{\bc(\ups) = \bc(\ups^*)\}$, $\ups(t) \stackrel{t \to \infty}{\longrightarrow} \ups^*$ exponentially, and (ii)  $\ups(t) \stackrel{t \to \infty}{\longrightarrow} \cA$ for other initial conditions. Further, if $\ups^* = \mathbf{0} \in \cR$, it is called \underline{q-attractor} if the above happens with  $ {\mathbb S}(\ups) := \{\pc = 0\}$.
\end{definition}
For systems modelling the BPs,  the following subset of the domain is relevant:
\begin{align}\label{invariant_set}
    \cD_I &:= \{\ups \in (\mathbb{R}^+)^4: \tc \leq \pc \leq \pa \mbox{ and } \ta \leq \pa\}.
\end{align}
It is not difficult to see that $\cD_I$ is positive invariant for ODE \eqref{eqn_ODE} - we prove the same in Theorem \ref{thrm_attractors_beta}.
Consider the following subset of $\cD_I$, which represents the combined domain of attraction towards $\cA\cup\cR$ (attractors and saddle points):
\begin{align}\label{eqn_domain of attraction}
    \cD &:= (\cD_\cA \cup \cD_\cR) \cap \cD_I  = \{\ups \in \cD_I: \ups(t) \to {\cA}\cup \cR \mbox{ as } t \to \infty, \mbox{ if } \ups(0) = \ups\}.
\end{align}
Thus, if the ODE starts in $\cD$, it converges asymptotically to $\cA\cup\cR$. The main result is:  when BP ($\Ups_n$) visits some compact subset of $\cD$ i.o., then either $\Ups_n$ converges asymptotically to $\cA\cup\cR$ or hovers around $\cR$ (notion defined below).
\begin{definition}\label{defn_hovers}
    The stochastic process $\Ups_n$ is said to  \underline{hover around a set} $\cR$ if $ \Ups_n \in N_\delta (\cR) \mbox{ i.o., for all } \delta >0 \mbox{ and  } \Ups_n \notin N_{\delta_1} (\cR) \mbox{ i.o., for some } \delta_1 >0$.
\end{definition}
\textit{Hovering around depicts a type of the limiting behavior of the stochastic process where the trajectory goes arbitrarily close to the set $\cR$ i.o., but still comes out of a neighbourhood of it i.o.}
Contrary to the existing results, our SA based Theorem \ref{thrm1} given below proves the possibility of above behavior as well as convergence to the saddle set ($\cR$). We require an extra assumption and the proof is deferred to the next section.
\begin{enumerate}[label=\textbf{A.\arabic*}, ref=\textbf{A.\arabic*}]
\setcounter{enumi}{3}
    \item Let $\cA\cap\cD_I$ be the attractor set  as in Definition \ref{defn_attractors}. Let each $\ups \in \cR\cap\cD_I$ be the q-attractor as in Definition \ref{defn_q_as}.  Consider $\cD$ as in \eqref{eqn_domain of attraction} and let  $\cS := \cD \cap \{\pa \leq b\}$, for some $b > 0$, be a compact subset of combined domain of attraction. Assume $p_{b} := P(\mathcal{V}) > 0$, where $\mathcal{V} :=  \{\omega : \Ups_n(\omega) \in \cS \mbox{ i.o.}\}$. \label{a4}
\end{enumerate}

\begin{theorem}\label{thrm1}
  Under \ref{a1}-\ref{a3}, we have:
  \begin{enumerate}[label=(\roman*)]
        \item For every $T>0$, a.s.  there exists a sub-sequence $(n_l)$ such that:
            $$
            \sup_{k: t_k \in [t_{n_l}, t_{n_l} + T]} d(\Ups_k, \ups(t_k - t_{n_l})) \to 0  \mbox{ as } l \to \infty, \mbox{ where } t_n := \sum_{k=1}^n \frac{1}{k} \mbox{ and}
            $$
        $\ups(\cdot)$ is the extended solution of ODE \eqref{eqn_ODE} which starts at $\ups(0) =
        \lim_{n_l \to \infty} \Ups_{n_l}$.
        \item Further, assume \ref{a4}. Then, $P({\cal C}_1 \cup {\cal C}_2) \geq p_b$, where
        \begin{align*}
                    \hspace{-10mm}
        \begin{aligned}
                 \hspace{7mm}{\cal C}_1 &: =\{\Ups_n \to (\cA \cup \cR)\cap \cD_I \mbox{ as } n \to \infty\}, \mbox{ and }
                {\cal C}_2 := \{ \Ups_n \mbox{ hovers around } \cR \}. \hspace{5mm}  \mbox{ \eop}
            \end{aligned}
        \end{align*}
  \end{enumerate}
\end{theorem}
Thus, the BP either converges to attractor/saddle set or it can hover around a saddle point, with probability at least $p_b > 0$; in fact, the saddle points are q-attractors defined in Definition \ref{defn_q_as}.  We will show that \ref{a1}-\ref{a4} are satisfied for BP with attack in Section \ref{sec_BPA} and proportion-dependent BP in Section \ref{sec_prop_BP}, with $p_b = 1$, i.e, the above results are true a.s.

\subsection{Significance of Theorem \ref{thrm1}} \label{subsec_significance_thrm1}

\noindent \textbf{BP trajectories -} Theorem \ref{thrm1}(i) provides \textit{a novel approach for studying the asymptotic trajectory of the BPs using ODE solution}.
Consider the solution of ODE \eqref{eqn_ODE} initialised with $\lim_{n_l \to \infty} \Ups_{n_l}$. Then, the BP $\Ups_k$ is  close to ODE solution $\ups(t_k-t_{n_l})$ at all transition epochs, $k$ with $t_k \in [t_{n_l}, t_{n_l} + T]$. This approximation improves as $n_l$ increases. \textit{The result is true a.s., for all $T < \infty$, independent of $p_b$ and only requires \ref{a1}-\ref{a3}.}  

We suggest a better finite-time approximation using a non-autonomous ODE in Section \ref{sec_6}, inspired by \cite{agarwal2022saturated}, where saturated  total population-dependent BP is studied.

\vspace{1mm}
\noindent \textbf{Limit proportion -}
Theorem \ref{thrm1}(ii) provides an alternate approach to derive limit behaviour via the attractors or saddle points (q-attractors) of ODE \eqref{eqn_ODE}.

In \textit{extinction paths}, where both populations get extinct, $\Ups_n \to \mathbf{0}$  as $n \to \infty$, say with probability $p_e > 0$. Thus, extinction paths are in the set $\mathcal{V}$ of \ref{a4}. While in survival paths, the BP either converges or hovers around $\left(\cA \cup (\cR - \{\mathbf{0}\})\right)\cap \cD_I$, with probability at least $p_b - p_e$. 
As an example of convergence to saddle point, the vector $\mathbf{0}$ is a saddle point of ODE \eqref{eqn_ODE} (shown in the proof of Theorem \ref{thrm_attractors_beta}) and is also a limit of the BP in extinction paths.

\vspace{1mm}
\noindent \textbf{Population independent to population dependent BPs -}
Say one is interested in the limit analysis of a BP whose limit mean matrix (say) $M^\infty$ matches with that of a  BP (call it known-BP), say again $M^\infty$, which is previously analyzed in the literature, e.g., the analysis of population-independent BPs are usually known, and can be  candidates for the known-BPs. In such scenarios, our approach could be helpful by using the following steps: i) one needs to derive the analysis of the corresponding ODE \eqref{eqn_ODE} with the help of the known limits of the known-BP, and ii) show that the BP under investigation visits the domain of attraction i.o., using the fact that the known-BP visits the same i.o.
We illustrate an easy procedure to analyze the attractor/repeller set of the ODEs and  the corresponding domain of attraction (which is sufficiently large to ensure $p_b = 1$) for  many example BPs further in the paper; we believe one can follow similar procedure to analyze the ODEs corresponding to known-BPs, with the additional knowledge of the limits of the known-BPs.

\vspace{1mm}
\noindent \textbf{Limitation -} By Theorem \ref{thrm1}, one can not comment on the individual probability of $\Ups_n$ converging to a particular limit in $\left(\cA \cup \cR\right)\cap \cD_I$ or the likelihood of hovering around. Further, $\Bc_n \to \{0, 1\}$ does not always imply the extinction of $x$ or $y$-type population. However, in BP with attack, this is true (see the discussion at the end of Appendix \ref{appendix_prelim}).

\section{Proof of Theorem \ref{thrm1}}\label{sec_proof_thrm1}

From equation \eqref{evolve_x_up_time}, the embedded process immediately after $n$-th death, when the death is for example of an $x$-type individual, is given by:
\begin{equation}\label{evolve_x_up_gen}
\begin{aligned}
\Cx_{n} &= \Cx_{n-1}  + \offs_{xx, n}(\Om_{n-1}) - 1, \ \ \  \Ax_{n} = \Ax_{n-1}  + \offs_{xx, n}(\Om_{n-1}), \\
\Cy_{n} &= \Cy_{n-1} + \offs_{xy, n}(\Om_{n-1}),  \ \ \ \Ay_{n} = \Ay_{n-1} + \offs_{xy, n}(\Om_{n-1}).
\end{aligned}
\end{equation}
To begin with, we make an important observation to derive an appropriate SA-based scheme which represents the above dynamics and also to prove a boundedness assumption for ratios $\Ups_n$ required for most SA-based studies. 

\noindent \textbf{Key idea:}
Consider a BP with population-independent and positive offspring, i.e., in \ref{a1}, assume $\offs_{ix}(\om) + \offs_{iy}(\om) = \overline{\offs}$ for all $\om$ and all $i \in \{x, y\}$. Let $\overline{\Pi}_n$ represent the sample mean formed by the sequence of   offspring plus the initial population size, i.e., 
\begin{align}\label{eqn_overline_S_n}
    \overline{\Pi}_n = \frac{1}{n}\left(\sum_{k=1}^n  \overline{\offs}_k + s_0^a \right).
\end{align}
By strong law of large numbers, $\overline{\Pi}_n \to \overline{m} := E[\overline{\offs}_1]$ a.s. 
For this special case, ${\Psi}^a_n = \overline{\Pi}_{n}1_{n < \nu_e} + \nu_e \overline{\Pi}_{\nu_e}/n 1_{n \geq \nu_e}$ (see \eqref{evolve_x_up_gen} and recall ${\Psi}^a_n = (\Ax_n+\Ay_n)/n$); hence ${\Psi}^a_n $ converges either to $0$ (in extinction paths, i.e., $\nu_e < \infty$) or to ${\overline m}$ (in survival paths); $\Pc_n$ respectively converges to $0$ or $\overline{m}-1$. This observation actually completes the  proof  for this special case with $\cA = \{(0, 0), (\overline{m}-1, \overline{m})\}$, further when single population (say $x$-type) is considered. It is well known that the sample mean \eqref{eqn_overline_S_n} \textit{can be written as a SA-based scheme} and in \eqref{eqn_SA_pi} given below, we will see that this is true even for the general case. Further, clearly, \eqref{eqn_overline_S_n} becomes an upper bound for all components of $\Ups_n$, which \textit{helps in bounding $\Ups_n$  uniformly in $n$ and a.s.} (see \eqref{Eqn_XnYnSnetc_general} given below), again under~\ref{a1}.  

Analogous to $\overline{\Pi}_n$ as in \eqref{eqn_overline_S_n}, one can construct a lower bounding sequence using $\underline{\offs}$ of \ref{a1}; this provides a uniform positive lower bound for $\Pc_n$, which will help the proof. 


{\bf Proof:}
For any $n \geq 1$, let $\Pi_n$ represent the sample mean formed by the sequence of (possibly $\om$-dependent) offspring plus the initial population size till $\nu_e$, i.e., 
\begin{align}\label{eqn_pi_n}
     \Pi_n &= \frac{1}{n}\left(\sum_{k=1}^{\min\{n, \nu_e\}} \left( H_k\offs_{x, k}(\Om_{k-1}) + \overline{H}_k \offs_{y, k}(\Om_{k-1})  \right)1_{\Pc_k > 0} + s_0^a \right), \mbox{ where } \nonumber  \\
    \offs_{x, k} &:= \offs_{xx, k} + \offs_{xy, k}, \ \  
    \offs_{y, k} := \offs_{yy, k} + \offs_{yx, k},  
\end{align}and $H_k= 1-\overline{H}_k$ is the indicator that an $x$-type individual dies at $k$-th epoch. It is easy to observe that ${\Pi}_n$  can be re-written as (observe that $\nu_e$ also equals $\inf\{n: \Pc_n = 0\}$)
\begin{align}\label{eqn_SA_pi}
\Pi_n =  \Pi_{n-1} + \frac{1}{n} \bigg[ \left(H_n \offs_{x, n}(\Om_{n-1}) + \overline{H}_n \offs_{y, n}(\Om_{n-1}) \right)1_{\Pc_{n-1} > 0} - \Pi_{n-1} \bigg].
\end{align}

In fact, $\Pi_n = \Pa_n$ for all $n\geq 1$, and so the above iterative equation represents $\Pa_n$. Similarly, other ratios in $\Ups_n$ can be re-written as (with $\epsilon_{n-1} := 1/n$, see \eqref{evolve_x_up_gen}, \eqref{eqn_pi_n}, \eqref{eqn_SA_pi}): 
\begin{align}\label{eq_stoch_approx_scheme_general}
\begin{aligned}
\Ups_n &= \Ups_{n-1} + \epsilon_{n-1}\mathbf{L}_{n-1}, \mbox{ where } \mathbf{L}_{n-1} := (L_{n-1}^{\psi, c}, L_{n-1}^{\theta, c}, L_{n-1}^{\psi, a}, L_{n-1}^{\theta, a})^t, \mbox{ with}\\
L_{n-1}^{\psi, c} &:=  \left\{ H_n\bigg(\offs_{x, n}(\Om_{n-1})-1\bigg) + \overline{H}_n \bigg(\offs_{y, n}(\Om_{n-1})-1\bigg)\right\} 1_{\Pc_{n-1} > 0}   -  \Pc_{n-1}, \\
L_{n-1}^{\theta, c} &:= \left\{H_{n}\bigg(\offs_{xx, n}(\Om_{n-1}) - 1\bigg) + \overline{H}_{n} \offs_{yx, n}(\Om_{n-1})\right\}1_{\Pc_{n-1} > 0} - \Tc_{n-1},  \\
L_{n-1}^{\psi, a} &:=  \bigg\{ H_n \offs_{x, n}(\Om_{n-1}) + \overline{H}_n \offs_{y, n}(\Om_{n-1}) \bigg\}1_{\Pc_{n-1} > 0}  - \Pa_{n-1}, \mbox{ and}\\
L_{n-1}^{\theta, a} &:= \bigg\{H_{n}\offs_{xx, n}(\Om_{n-1}) + \overline{H}_{n} \offs_{yx, n}(\Om_{n-1})\bigg\}1_{\Pc_{n-1} > 0}  - \Ta_{n-1}. 
\end{aligned}
\end{align}


The proof of part (i) has two major steps: (a) to construct a sequence of piece-wise constant interpolated trajectories for almost all sample-paths; (b) to prove that the designed trajectories are equicontinuous in extended sense\footnote{\label{footnote_defn_equi}\begin{definition}
  {\bf Equicontinuous in extended sense (\cite[Equation (2.2), pp. 102]{kushner2003stochastic})):} \label{defn_equi} Suppose that for each $n$, $f_n(\cdot)$
is an $\mathbb{R}^r$-valued measurable function on $(-\infty,\infty)$ and $(f_n(0))$ is bounded.
Also suppose that for each $T$ and $\epsilon > 0$, there is a $\delta > 0$ such that
\begin{align}\label{eqn_footnote}
\limsup_n \sup_{0\leq t-s\leq \delta, |t| \leq T} |f_n(t) - f_n(s)| \leq \epsilon. 
\end{align}
 Then the sequence $(f_n(\cdot))$ is said to be equicontinuous in the extended sense.
\end{definition}}. These steps are majorly as in \cite[Theorems 2.1-2.2]{kushner2003stochastic}, but for the changes required for measurable $\mathbf{\ga}(\cdot)$.

For many steps of the proof, we will work only with $\tc$-component of the vector $\ups$, when the proof for the remaining components goes through in exactly similar manner.

Let $\Ups^n(\cdot) := (\Psi^{n, c}(\cdot), \Theta^{n, c}(\cdot), \Psi^{n, a}(\cdot), \Theta^{n, a}(\cdot))$  be  the constant piece-wise interpolated trajectory defined as below (see \eqref{eq_stoch_approx_scheme_general}, and recall $t_n = \sum_{i=1}^n \epsilon_{i-1}$):
\begin{align}\label{eqn_interpolated_traj_1}
    \Theta^{n, c}(t) := \Tc_n + \sum_{i=n}^{\eta(t_n + t)-1} \epsilon_i L_i^{\theta, c}, \mbox{ for all } t \ge 0,
\end{align}
$\Psi^{n, c}(t), \Psi^{n, a}(t)$ and $\Theta^{n, a}(t)$ are defined analogously. Towards proving equicontinuity, we first consider upper-boundedness of $\Ups^n(0) = \Ups_n$, as the iterates are trivially lower bounded by $0$. The claim is immediately true by strong law of large numbers a.s., to be more precise on the set $\{\overline{\Pi}_n  \to \overline{m}\}$,  because of the following  observation (see \eqref{eqn_overline_S_n}-\eqref{eq_stoch_approx_scheme_general}):
\begin{equation}\label{Eqn_XnYnSnetc_general}
     \Pc_n \leq \Pa_n \mbox{ and } \Tc_n \leq \Ta_n \le \Pa_n = \Pi_n \leq  \overline{\Pi}_n   \mbox{ for all } n, 
\end{equation}
as then for any sample path $\omega \in \{\overline{\Pi}_n  \to \overline{m}\}$ and   $\epsilon > 0$, there exists a $n_\epsilon(\omega) < \infty $,
\begin{align}\label{eqn_bounded_iterates}
\begin{aligned}
\hspace{6mm}
   \sup_n  \max\{\Theta^{n, c}(0), \Psi^{n, c}(0), \Theta^{n, a}(0), \Psi^{n, a}(0)\}  \leq \hspace{16mm} &    \\
   &\hspace{-8cm}\max\left \{ \max_{n \le n_\epsilon(\omega)\}}  \max\{\Theta^{n, c}(0), \Psi^{n, c}(0), \Theta^{n, a}(0), \Psi^{n, a}(0)\}, \   \overline{m} + \epsilon \right \}. 
   \end{aligned}
\end{align}
%
%
Towards the second part of equicontinuity (see \eqref{eqn_footnote} in footnote \ref{footnote_defn_equi}), the interpolated trajectory for $\Theta^{n,c}(\cdot)$ in \eqref{eqn_interpolated_traj_1} can be re-written in `almost integral form', for any $t \geq 0$:
\begin{align}\label{eqn_diff_term1}
\begin{aligned}
     \Theta^{n, c}(t) &:= \Tc_n + \int_0^t \rho_\theta^c(\Ups^n(s), s) ds + {\cal E}_1^{n, c}(t), \mbox{ with the difference term, } \\
     {\cal E}_1^{n, c}(t) &:= \sum_{i=n}^{\eta(t_n + t)-1} \epsilon_i L_i^{\theta, c} - \int_0^t \rho_\theta^c(\Ups^n(s), s) ds, \mbox{ where}
\end{aligned}
\end{align}
$\gna = (\rho_\psi^c, \rho_\theta^c, \rho_\psi^a, \rho_\theta^a)$ is the conditional expectation, $E[\mathbf{L}_n|\mathcal{F}_n] =: \gna(\Ups_n, t_n)$, with respect to the sigma algebra, ${\cal F}_n  := \sigma\{\Om_k : 1 \leq k < n \}$, and equals (see \eqref{eq_stoch_approx_scheme_general}):

 \vspace{-4mm}
{\small
\begin{align}\label{eqn_nonauto_g}
\rho_\psi^c (\ups, t) &:=      \left\{ \bc \bigg(m_{xx}(\om) +     m_{xy}(\om)\bigg) + (1-\bc)\bigg(m_{yy}(\om) + m_{yx}(\om)\bigg) - 1 \right\} 1_{\{\pc > 0\}}- \pc, \nonumber \\
\rho_\theta^c (\ups, t)  &: = \left\{ \bc \bigg(m_{xx}(\om) - 1\bigg) + (1-\bc) m_{yx}(\om)\right\} 1_{\{\pc > 0\}} - \tc,  \\
\rho_\psi^a (\ups, t) &:=      \left\{ \bc \bigg(m_{xx}(\om) + m_{xy}(\om)\bigg) + (1-\bc)\bigg(m_{yy}(\om) + m_{yx}(\om)\bigg)\right\} 1_{\{\pc > 0\}} - \pa,   \mbox{ \normalsize  and,}  \nonumber \\
\rho_\theta^a (\ups, t)  &: = \bigg\{\bc m_{xx}(\om) + (1-\bc) m_{yx}(\om)\bigg\} 1_{\{\pc > 0\}} - \ta,   \mbox{ \normalsize with, }   \eta(t) := \max\left  \{ n: t_n \le t \right \},  \nonumber\\
\om &= \om(\ups, t) := 
  \big(\tc \eta(t),  \ (\pc-\tc)\eta(t),\ \ta \eta(t), \  (\pa-\ta)\eta(t)\big). \nonumber
\end{align}}

We further re-write the interpolated trajectory using the autonomous ODE \eqref{eqn_ODE}:
\begin{align}\label{eqn_diff_term2}
\begin{aligned}
     \Theta^{n, c}(t) &:= \Tc_n + \int_0^t g_\theta^c(\Ups^n(s)) ds + {\cal E}_1^{n, c}(t) + {\cal E}_2^{n, c}(t), \mbox{ where}\\
     {\cal E}_2^{n, c}(t) &:= \int_0^t \rho_\theta^c(\Ups^n(s), s) ds - \int_0^t g_\theta^c(\Ups^n(s)) ds.
\end{aligned}
\end{align}
In Appendix \ref{appendix_B}, we show that ${\cal E}_1^{n, c}(t) + {\cal E}_2^{n, c}(t)$ converges uniformly to $0$, as $n \to \infty$, over any finite time window and further show:

\begin{lemma}\label{lemma_equi_cont_thrm1}
The sequence $(\Ups^n(\cdot))$ is equicontinuous in extended sense a.s. \eop
\end{lemma}
Now, consider the set $N$ of all sample paths for which $(\Ups^n(\cdot))$ is not equicontinuous - by Lemma \ref{lemma_equi_cont_thrm1}, $P(N) = 0$ (see proof of above Lemma for precise definition of $N$). Then, by extended version of Arzela-Ascoli Theorem \cite[section 4, Theorem 2.2, pp. 127]{kushner2003stochastic}, there exists a sub-sequence $(\Ups^{n_m}(\omega, \cdot))$ which converges to some continuous limit, call $\ups(\omega, \cdot)$, uniformly on each bounded interval for $\omega \notin N$ such that: 
\begin{align}\label{eqn_ups_inf}
\ups(t) = \lim_{n_m \to \infty} \ups_{n_m}(\omega) + \int_0^t \ga (\ups(s)) ds.
\end{align}
Thus, for every $\epsilon > 0$ and $T > 0$, there exists $n(\omega, \epsilon, T)$ such that: 
\begin{align}\label{eqn_dist_scheme_ODE_}
\sup_{l \in L} d({\Ups}_l , \ups(t_l - t_{n_m})) \leq \epsilon/2 \mbox{ for all } n_m \geq n(\omega, \epsilon, T), 
\end{align}
where $L:= \{l :  t_{n_m} \leq t_l \leq T + t_{n_m}\}$; observe for any $l \in L$, ${\Ups}^{n_m}(t) = {\Ups}_l$ if $t = t_l - t_{n_m}$. Now, we are left to show that $\ups(\cdot)$ in \eqref{eqn_ups_inf}, the solution of the fixed point equation (of the integral operator), is the extended solution of ODE \eqref{eqn_ODE} starting at $\ups(0) = \lim_{n_l \to \infty} \Ups_{n_l}$, i.e., 
$$
    \lim_{h \to 0} \frac{\ups(t+h) - \ups(t)}{h} = \ga(\ups(t)) = \frac{d \ups(t)}{d t} \mbox{ for almost all } t.
$$
One can easily show that the function $\ga \circ \ups$ is locally integrable, and thus, by 
\cite[Theorem 3.21]{folland1999real}, the claim holds. 
This completes part (i).

For part (ii), under \ref{a4}, the proof is again inspired from \cite{kushner2003stochastic} and \cite[Theorem 2.3.1, pp. 39]{kushner2012stochastic}, even when the solution of ODE \eqref{eqn_ODE} is in extended sense, not the classical one. Further major difference in the proof is to include the arguments required to prove the event of hovering around $\cR$. We complete this proof in Appendix \ref{proof_thrm1}. \eop

\section{\textbf{Derivation of $\cA$ and $\cR$ - analysis of proportion ODE}}\label{procedureAR}

Under \ref{a2}, $\om$-dependent mean functions converge to just $\bc$-dependent mean functions, and thus, one may anticipate that the analysis of $ \bc(\ups(t))  = \bc(t)$ plays a crucial role. In fact, we claim and prove that the time limits of $\bc$, obtained from the following  limit ODE for $\bc$ (derived using \eqref{eqn_ODE}), leads to the required analysis:

\vspace{-6mm}
\begin{align}\label{eqn_beta_ODE}
\begin{aligned}
\dot{\bc} &= \frac{1}{\pc} g_\beta(\bc)1_{\{\pc > 0\}},\mbox{ where}\\ 
g_\beta(\bc) &:= - \bc \minf_{xy}(\bc) + (1-\bc) \minf_{yx}(\bc)\\
    &\hspace{15mm}+ \bc(1-\bc) \bigg\{\minf_{xx}(\bc) + \minf_{xy}(\bc) - \big(\minf_{yx}(\bc) +\minf_{yy}(\bc)\big)  \bigg\}.
    \end{aligned}
\end{align}
From above, $g_\beta$ depends only on $\bc$, thus, \textit{one might expect that the asymptotic analysis of $\bc$ is independent of other components of $\ups$}. We will see that this is indeed true, and in fact, asymptotic analysis of all components of $\ups$ can be derived using $g_\beta$.
In this regard, we define the following:
\begin{definition}
    Any point $\bstar \in [0,1]$ is  (projected) \underline{p-stable} if $\mathbf{h}(\bstar)$ is an attractor for ODE \eqref{eqn_ODE}; a $\bstar$ is called \underline{p-saddle} if $\mathbf{h}(\bstar)$ is a saddle point, more specifically, q-attractor.
\end{definition}

Under certain conditions, we will show that the attractors of the following one-dimensional ODE are p-stable, while the repellers\footnote{\label{defn_repeller}Any point $\bstar \in [0,1]$ is called a repeller of ODE \eqref{eqn_beta_ode_simple} if $g_\beta(\bstar) = 0$ and $\bc(t) \nto \bstar$ as $t \to \infty$ when $\bc(0) \in {\cal N}_\epsilon(\bstar)$ for some $\epsilon > 0$.} are p-saddle:
\begin{align}\label{eqn_beta_ode_simple}
    \dot{\bc} &= g_\beta(\bc).
\end{align}

\begin{wrapfigure}{R}{0.36\textwidth}
  \vspace{-5mm}
    \centering
    \includegraphics[trim = {2.8cm 0.1cm 3cm 0cm}, clip, scale = 0.24]{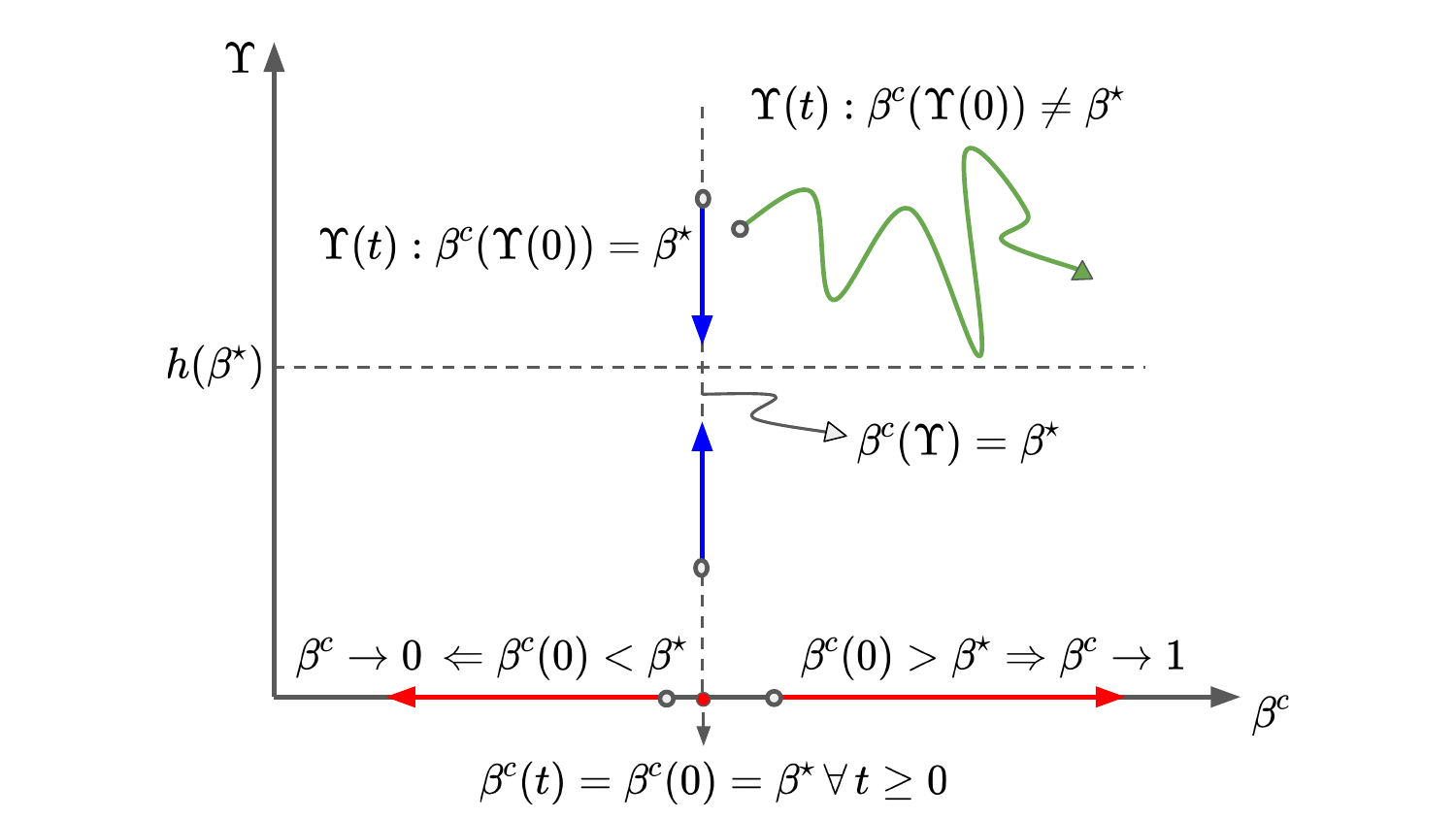}
    \vspace{-10mm}
    \caption{Repeller of \eqref{eqn_beta_ode_simple} leads to saddle point of \eqref{eqn_ODE}}\label{fig_saddle}
    \vspace{-6mm}
\end{wrapfigure}
When $\bstar$ is a repeller of \eqref{eqn_beta_ode_simple}, we have $g_\beta(\bstar) = 0$. Thus, when ODE \eqref{eqn_ODE} is initialised with $\bc(\ups(0)) = \bstar$, the ODE solution may remain in affine sub-space $\{\bc(\ups) = \bstar\}$ and may converge to $\mathbf{h}(\bstar)$ (see figure \ref{fig_saddle}). 
While if $\bc(\ups(0)) \neq \bstar$,  one might expect the solution of ODE \eqref{eqn_ODE} to repel away from $\mathbf{h}(\bstar)$, by definition of repeller. These observations indicate that $\bstar$ should be p-saddle and we precisely prove the same in our \textit{second important result} below. This result is instrumental in deriving $\cA$ and $\cR$ using the limit set of ODE \eqref{eqn_beta_ode_simple}; see Appendix \ref{proof_thrm2} for the proof. 

\begin{theorem}\label{thrm_attractors_beta}
Let ${\cal I} = \{x_i^* : 1 \leq i \leq n\}$ be the set of dis-continuities with $1 \le n < \infty$ and ${\cal J} :=  \{y_i^* : 1 \leq i \leq m\}  \subset {\cal I}^\complement$ be the set of points with $m < \infty$  (${\cal J}$ is empty when $m = 0$) such that:

    \noindent  (a) $g_\beta(x) = 0$ for each $x \in {\cal I} \cup {\cal J}$, i.e., ${\cal I} \cup {\cal J}$ is the set of equilibrium points for \eqref{eqn_beta_ode_simple},
    
    \noindent (b) for each $1\leq i \leq n$, there exists an open/closed/half-open non-empty interval around $x_i^* \in {\cal I}$, say ${\cal N}_i^*$, such that 
    
    (i) $\cup_{1\leq i\leq n} {\cal N}_i^* = [0,1]- {\cal J}$ and  ${\cal N}_i^* \cap {\cal N}_j^* = \emptyset$ for $i\neq j$,
    
    (ii) $g_\beta(\beta) > 0$ for all $\beta \in  {\cal N}^{-}_i:= {\cal N}^*_i\cap[0, x_i^*)$, $g_\beta$ is Lipschitz continuous on $ {\cal N}^{-}_i$,

    (iii) $g_\beta(\beta) < 0$ for all $\beta \in {\cal N}^{+}_i:
    = {\cal N}^*_i \cap (x_i^*, 1]$, $g_\beta$ is Lipschitz continuous on ${\cal N}^{+}_i$.

\noindent  Then, ODE \eqref{eqn_ODE} satisfies \ref{a3}. Further, the set ${\cal I}$ is an attractor for \eqref{eqn_beta_ode_simple} and p-stable for  \eqref{eqn_ODE}; also, ${\cal J}$ is the set of repellers for \eqref{eqn_beta_ode_simple} and p-saddle for \eqref{eqn_ODE}. 
Furthermore,  $\cA := \{\mathbf{h}(x_i^*) : x_i^* \in {\cal I}\}$ is the attractor set, $\cR := \{\mathbf{h}(y_i^*) : y_i^* \in {\cal J}\}\cup \{\mathbf{0}\}$ is the saddle set  in $\cD_I$ and entire $\cD_I$ is the combined domain of attraction for  \eqref{eqn_ODE}. \eop
\end{theorem}
We believe that the above Theorem can be extended for $g_\beta$ which is continuous, by standard ODE results, and we leave it for future. However, we require $g_\beta$ to be discontinuous for BP with attack (see assumption \ref{k2} in Section \ref{sec_BPA}), and thus the hypothesis of Theorem \ref{thrm_attractors_beta}. The last part of the Theorem asserts that the p-stable/p-saddle points are the only attractors/saddle points of ODE \eqref{eqn_ODE}, other than $\mathbf{0} \in \cR$.

\section{Related work and extensions}\label{sec_survey}
There is a vast literature related to BPs, however, we simply discuss few relevant strands related to our work. In sub-section \ref{desc_specialCases}, we discuss a variety of existing BPs and their extensions to illustrate the generality of our result. 

Irreducible population-dependent BP with discrete and continuous-time framework are considered in \cite{klebaner1989geometric, jagers1997coupling} respectively; they do not consider total population-dependent offspring; further, the population-dependent mean matrix converges to a constant mean matrix, but we support proportion-dependent mean matrix in the limit.  In \cite{athreya1968some}, authors consider continuous-time, but population-independent, irreducible BPs. 

In \cite{coffey1991galton}, the prey-predator BP is analyzed in discrete-time setting and co-survival conditions are identified, but the limit proportion is not derived; they also do not consider population-dependency.  In Section \ref{sec_BPA}, we consider a continuous-time population-dependent BP with double-sided attack and acquisition. One can also analyze the continuous-time population-dependent variant of prey-predator BP using our results, as illustrated below in sub-section \ref{desc_specialCases}.

In \cite{agarwal2021co}, we introduce BP with attack and provide limit proportion for the case with population-independent and symmetric offspring, i.e., with $m_{xx}(\om) = m_{yy}(\om) = m$, for all $\om$, for some $m > 1$. We significantly generalize by considering total population-dependency and symmetric/asymmetric offspring. We analyze a particular case of proportion-dependent BP (offspring depend on the proportion of the populations) along with other co-authors in \cite{kapsikar2020controlling}. Our results cover the model in \cite{kapsikar2020controlling} and can also be used to generalize their result which will be a part of our future work.

\textbf{\polya urn models:} In \cite{athreya1968embedding}, it is shown that the \polya urn models can be embedded into a continuous-time population-independent BP. Thus, the asymptotic analysis of the continuous-time BPs can be derived using the corresponding analysis of the \polya urn models. However, our work differs from the \polya urn literature (\cite{higueras2006central, janson2004functional, arthur1987non}) for reasons mentioned in the introduction - they neither consider total population-dependency nor commonly deal with extinction (non-replacement) scenarios; recall, in BPs, extinction occurs with non-zero probability, even in the super-critical regime. In \cite{janson2004functional}, which is an exception, the possibility of extinction is considered, but they do not consider population-dependency.

In \cite{higueras2006central}, authors analyze the urn model with the removal of balls of other colours (not the chosen one) - same as a negative offspring in BP with attack of Section \ref{sec_BPA}. However, they assume a unique attractor for ODE and a constant number of additions (offspring) to the urn. We again have a significant generalization with a random number of offspring and where the random trajectory of the BP with attack can converge 
 to or hover around one of the attractors/saddle points of ODE (see Corollary \ref{corollary_BPA}).

\subsection{Existing BPs} \label{desc_specialCases}
We now discuss a wide range of existing and new BPs captured via dynamics provided in \eqref{evolve_x_up_time}. By Theorem \ref{thrm1}, the reader can observe that any BP, which is either population dependent/independent, or irreducible/decomposable, or have positive/negative offsprings, or have any combination of these and more during transience, can be analyzed if one analyzes the limiting ODE \eqref{eqn_ODE}. In fact, the analysis is governed by the ODE \eqref{eqn_beta_ode_simple}  constructed only using limit mean matrix $M^\infty(\bc)$, by following the steps given in section \ref{procedureAR}. This procedure is validated by Theorems \ref{thrmBPA} and \ref{thrmProp}  provided in the coming sections, which also pave the path for showing $p_b = 1$, i.e., almost sure convergence in Theorem \ref{thrm1}. Recall, we would also require all such  processes to satisfy bounding and convergence rate assumptions as in \ref{a1} and \ref{a2} respectively.


 \subsubsection{\textbf{Prey-predator BP}}\label{preyBP}   As said before, the authors in \cite{alsmeyer1993galton, coffey1991galton} study a population-independent prey-predator BP in discrete-time scenario; obviously the attacks can not exceed the population available to be attacked. With $\bpam_{xx}$, $\bpam_{yy} > 1$ representing the population-independent mean offsprings of predator and prey respectively, it is concluded in \cite{coffey1991galton} that:
    
    (i) if  $\bpam_{xx} \geq \bpam_{yy}$, then prey gets extinct w.p. $1$, if predator survives, 
    
    (ii) otherwise, for any $\cx(0)$, there exists $\cy(0)$ such that $P_0(\Cx \nto 0, \Cy \nto 0) > 0$. 

\noindent To exactly achieve the continuous-time version of the prey-predator BP, for each $\Om$, let $\offs_{ii}(\Om) = \xi_{ii} \geq 0$, for every $i\in \{x, y\}$, $\offs_{xy}(\Om) = -\xi_{xy}(C^y)  \leq 0$ and $\offs_{yx}(\Om) = 0$ in  \eqref{evolve_x_up_time}, which further satisfy \ref{a1}. We discuss the derivation of limit mean matrix for BPs with negative offsprings in Appendix \ref{appendix_prelim}, and similarly consider: 
    
    \vspace{-4mm}
    {\small
    \begin{align}\label{mean_prey}
    M(\bc) := 
    \begin{bmatrix}
    \bpam_{xx} & -\bpam_{xy} 1_{\{\bc < 1\}}\\
    0 & \bpam_{yy}
    \end{bmatrix}, \mbox{ where } \bpam_{ii} = E[\xi_{ii}], \mbox{ and } 
    \end{align}}the expected number of attacks, $E[\offs_{xy}(\Om)] = -\bpam_{xy}1_{\{\bc < 1\}}$. Now, under the procedure given in section \ref{procedureAR} and Theorem \ref{thrm1}, we get that\footnote{the only difference is that in Theorem \ref{thrmBPA}, $\offs_{ii}(\Om) = \xi_{ii}(\Om) + \xi_{ij}(C^j)$ (includes acquisition), in comparison to $\offs_{ii}(\Om) = \xi_{ii}$. For prey-predator BP, following the steps as in proof of Theorem \ref{thrmBPA}, one can get that if $\bpam_{xx}^\infty < \bpam_{yy}^\infty$, then $\cA = \{\mathbf{h}(0), \mathbf{h}(1)\}$ and $\cR = \{\mathbf{h}(\bstar_r), \mathbf{0}\}$; otherwise, $\cA = \{\mathbf{h}(1)\}$ and $\cR = \{\mathbf{h}(0), \mathbf{0}\}$. Hence, $p_b = 1$: it is not difficult to show that the dynamics \eqref{evolve_x_up_time} visits the respective DoA i.o., a.s.}: 
    
    (i) $\Bc_n \to \{0, 1, \bstar_r\}$ 
 or hovers around $\bstar_r$ a.s. for some $\bstar_r \in (0, 1)$ such that $\ga_\beta(\bstar_r) = 0$, if $\bpam_{xx} < \bpam_{yy}$, and
    
    (ii) $\Bc_n \to \{0, 1\}$ or hovers around $0$ a.s. if $\bpam_{xx} \geq \bpam_{yy}$. 
    
    \noindent These observations mirror the conclusions in \cite{coffey1991galton} stated above, except for the hovering around aspect; we additionally provide the limit proportion of the two populations ($\bstar_r$) in the co-survival sample paths. Furthermore, \textit{these results are true even for continuous-time population-dependent variant of prey-predator BP, if at the limit $M^\infty(\bc) = M(\bc)$ given in \eqref{mean_prey}.} It should now be clear that using our approach, more complex BPs can be analyzed, e.g., prey-predator population-dependent BPs with self offsprings also being proportion-dependent (i.e., $e_{ii}(\bc)$).

\subsubsection{\textbf{Irreducible and decomposable BPs}}  
Consider the population-dependent BPs (with positive or negative offsprings in transience) which converge to  irreducible or decomposable proportion dependent or independent BPs.  The dynamics of such BPs can be modelled by letting $\offs_{ij}(\om) \approx \offs_{ij}(\bc)$ for each $i, j$ `for large values of $\om$' in \eqref{evolve_x_up_time}. Now, we discuss the analysis of such BPs:

    (i) Consider the population-dependent BP with constant (decomposable/irreducible) limit mean matrix. Further, consider a population-independent BP with the mean matrix as the said limit mean matrix throughout. Many such population-independent variants are analyzed in the literature, which can assist in the analysis of the population-dependent BP, as in sub-section \ref{subsec_significance_thrm1}.

    (ii) The classical approach uses the martingales to analyze the population-independent BPs, where appropriate rates need to be deduced; there is no formal procedure to compute such rates for decomposable BPs. However, if one adapts to our method as explained in section \ref{procedureAR}, then analysis of the ODE \eqref{eqn_beta_ode_simple} suffices.
    
     (iii)  As an example, consider a decomposable population-independent variant of BP where only $x$-type parent can produce offsprings of both types such that the mean matrix is given by: 
     \begin{align}\label{mean_irr_dec}
     M =
     \begin{bmatrix}
     m_{xx} & m_{xy}\\
     0 & m_{yy}
     \end{bmatrix}, \mbox{ where }  m_{xy} > 0.
     \end{align}
     Let $1 < m_{xx} < m_{yy}$. Then, according to \cite[Theorem 1(ii)]{ranbir2019decomposable}, \cite{kesten1967limit} for $W_1, W_2 \geq 0$:
    $$
    \Cx(t) e^{-m_{xx} t} \to W_1 \mbox{ a.s., and }\Cy(t) e^{-m_{yy} t} - \frac{m_{xy}}{m_{yy} - m_{xx}}\Cx(t) e^{-m_{yy} t}   \to W_2 \mbox{ a.s.}
    $$
    Thus, $\Cx$, $\Cy$ grow at rates $m_{xx}$, $m_{yy}$ respectively; hence, proportion $\Bc(t) \to 0$ as $t \to \infty$. The same result can be derived using our approach, from the solution of the ODE \eqref{eqn_beta_ode_simple} which simplifies to the following for the example in \eqref{mean_irr_dec}:
    $$
    \dot{\bc} = \bc(1 - \bc)\left( m_{xx} - m_{yy} \right) - (\bc)^2 m_{xy}.
    $$
    The solution of the above ODE is:
    $$
    \bc(t) = \frac{ak}{kb - e^{-at}}, \mbox{ where } a = m_{xx} - m_{yy}, \ b = a + m_{xy} \mbox{ and } k = \frac{\bc(0)}{b\bc(0) - a}.
    $$
    Hence, clearly $\bc(t) \to \bstar =  0$, as $t\to  \infty$ (for any initial condition), and so, $\Bc_k = \Bc(t_k) \to 0$, i.e., $\Ups_k \to \{\mathbf{h}(0), \mathbf{0}\}$, or $\Ups_k$ hovers around $\{\mathbf{0}\}$ as $k\to  \infty$ a.s.\footnote{For decomposable population-independent BP, following the steps as in proof of Theorem \ref{thrmProp}, one can get that $\cA = \{\mathbf{h}(0)\}$, and repeller set, $\cR = \{\mathbf{0}\}$. Hence, $p_b = 1$: it is not difficult to show that the dynamics \eqref{evolve_x_up_time} visits the respective domain of attraction i.o.} 
    
    Next let $m_{xx} > m_{yy} > 1$. Consider the Markovian variant of the age-dependent decomposable population-independent two-type process in \cite{jagers1969proportions}. Then, the authors prove that the limit proportion is given by:
    \begin{align}\label{eqn_beta_decomposable_pnd}
    \lim_{t \to \infty} \Bc(t) = \frac{m_{xx} - m_{yy}}{m_{xx} - m_{yy} + m_{xy}} \mbox{ a.s.}
    \end{align}
    Using our analysis, one can derive the same result for the continuous-time and Markovian variant, see \eqref{eqn_beta_decomposable}. We re-iterate, one can derive the analysis of population-dependent BPs whose limit mean matrix  has the structure as in \eqref{mean_irr_dec}, as described in sub-section \ref{subsec_significance_thrm1}. In fact, we provide a far more general result in sub-section \ref{sec_prop_BP}, and the precise details of the model are discussed in section \ref{sec_prop_BP}.

\section{Branching Process with Attack} \label{sec_BPA}
Consider a BP with two population types, say $x$ and $y$. Each individual of any type lives for a random time, $\tau \sim exp(\lambda)$, where $\lambda \in (0, \infty)$. It produces a random number of offspring before dying. The BP also includes  attack and acquisition by rival types. 
    
To be precise, an individual of (say) $x$-type produces $\xi_{xx}(\Om(\tau^-))$ offspring of its type. Further, it attacks/removes $\xi_{xy}(C^y(\tau^-))$ individuals of $y$-type population; naturally, the attacked population can not exceed the population available to be attacked at $\tau^-$, hence $\xi_{xy}(C^y(\tau^-)) \leq C^y(\tau^-)$ a.s.; note that the number of attacks do not depend on the size of the attacking population. The attacked individuals are then deleted from the $y$-population, and acquired by (i.e., added to) the $x$-population.  Thus, for example, when a $x$-type individual dies, the current populations change as follows:
    
    \vspace{-0.5cm}
    {\small
    \begin{align*}
    \Cx(\tau^+) &= \Cx(\tau^-) + \xi_{xx}(\Om(\tau^-)) + \xi_{xy}(C^y(\tau^-)) - 1,  \mbox{ and }
    \Cy(\tau^+) = \Cy(\tau^-) - \xi_{xy}(C^y(\tau^-)).
    \end{align*}}The total and $y$-population also evolve similarly. 
    We call such a BP as \textit{Branching Process with Attack}.  The dynamics in  \eqref{evolve_x_up_time} capture this BP, when for each $i, j$:
    \begin{align}\label{eqn_dynamics_BPA}
    \offs_{ii}(\Om(\tau^-)) := \xi_{ii}(\Om(\tau^-)) + \xi_{ij}(C^j(\tau^-)), \mbox{ \normalsize and } \offs_{ij}(\Om(\tau^-)) := - \xi_{ij}(C^j(\tau^-)).
    \end{align}
Next, we assume:
\begin{enumerate}[label=\textbf{K.\arabic*},ref=\textbf{K.\arabic*}]
    \item For each $i \in \{x, y\}$, assume that there exist integrable random variables, $\overline{\xi}$, $\underline{\xi}$, such that $0 \leq \underline{\xi} \leq \xi_{ii}(\om) \leq \overline{\xi}$ a.s. for each $\om$ and $E[\overline{\xi}]^2 < \infty$, $E[\underline{\xi}] > 1$. Further, let the attack offspring $\xi_{ij}(\om)$ be integrable for each $\om$ and for each $i \neq j \in \{x, y\}$. 
    \label{k1}
\end{enumerate}
    The above assumption immediately implies \ref{a1}. Define the expectations conditioned on $\om$ as $\bpam_{ij}(\om) := E[\xi_{ij}(\om)]$ for  $i, j \in \{x, y\}$. We further assume (see \eqref{eqn_dynamics_BPA}):
\begin{enumerate}[label=\textbf{K.\arabic*},ref=\textbf{K.\arabic*}]
\setcounter{enumi}{1}
    \item For  $i, j \in \{x, y\}$, let $\bpam_{ij}^\infty \geq 0$  with $\bpam_{xy}^\infty > 0$. Assume $\minf_{ij}(\bc)$ satisfy the following: \label{k2}
    \begin{align*}
    \begin{aligned}
    \minf_{xy}(\bc) &=  -\bpam_{xy}^\infty 1_{\{\bc < 1\}}, \ \minf_{yx}(\bc) = -\bpam_{yx}^\infty  1_{\{\bc > 0\}},\\
    \minf_{xx}(\bc) &= \bpam_{xx}^\infty  + \bpam_{xy}^\infty  1_{\{\bc < 1\}} \mbox{ and }  \minf_{yy}(\bc) = \bpam_{yy}^\infty  + \bpam_{yx}^\infty  1_{\{\bc > 0\}}.
    \end{aligned}
    \end{align*}
    Further, assume the conditions of \ref{a2} are satisfied with $\{(m_{ij}, m_{ij}^\infty)\}_{i, j}$ replaced by $\{(\bpam_{ij}, \bpam_{ij}^\infty)\}_{i, j}$.
\end{enumerate}
We are interested in the BP where attack is prominent\footnote{If both $\bpam_{xy}^\infty, \bpam_{yx}^\infty = 0$, then it will lead to two independent (non-attacking) BPs at limit; if required, one can derive the analysis for this case, as done in Theorem \ref{thrmBPA}.} even at the limit, thus, $\bpam_{xy}^\infty > 0$ without loss of generality in \ref{k2}. If $\bpam_{yx}^\infty = 0$, then it leads to single-sided attack at limit, but recall anything is possible in transience. 
Observe the cross-mean function in \ref{k1} converge to (almost) constant limit, e.g., $\bpam_{xy}(\om) \stackrel{s^c \to \infty}{\to} \bpam_{xy}^\infty1_{\{\bc < 1\}}$. The reason behind the indicator is that there is no attack at limit when $\bc = 1$; this is because $\Cy_n \to 0$ when $\limsup_{n \to \infty} \bc(\Ups_n) = 1$ as proved at the end of Appendix \ref{appendix_prelim}.

For BP with attack, the ODE \eqref{eqn_ODE} has the following form: 

\vspace{-4mm}
{\small
\begin{align}\label{ODE_BPA}
    \begin{aligned}
        \dot{\ups} &= \mathbf{h}(\bc) 1_{\{\pc > 0\}} - \ups, \mbox{ where } \mathbf{h}(\bc) := (h_\psi^c, h_\theta^c, h_\psi^a, h_\theta^a), \mbox{ is such that}\\
        h_\psi^c &=  \bc \bpam_{xx}^\infty + \big(1 -\bc \big) \bpam_{yy}^\infty - 1, \
        h_\theta^c = \bc\left(\bpam_{xx}^\infty +\bpam_{xy}^\infty1_{\{\bc < 1\}} -1\right) - \left(1-\bc\right)\bpam_{yx}^\infty 1_{\{\bc > 0\}}, \\
        h_\psi^a &=  \bc \bpam_{xx}^\infty + \big(1 -\bc \big) \bpam_{yy}^\infty, \mbox{ and } h_\theta^a = \bc\left(\bpam_{xx}^\infty+\bpam_{xy}^\infty1_{\{\bc < 1\}}   \right) - \left(1-\bc\right)\bpam_{yx}^\infty 1_{\{\bc > 0\}}.
    \end{aligned}
\end{align}}

 We begin with the analysis of the above ODE towards providing ODE approximation result for BP with attack using Theorem \ref{thrm1}.
 
\subsection{Analysis of ODE for BP with attack} 
Define the parameter vector $\mathbf{e} := \{\bpam_{ij}^\infty: i, j \in \{x, y\}\}$, and consider the following class of limit mean functions (by \ref{k2}, the vector ${\mathbf e}$ defines $M^\infty$):
\begin{align}\label{eqn_setE}
    \cM &:= \{ \mathbf{e}: \bpam_{yx}^\infty  > 0 \} \cup \{ \mathbf{e}:  \bpam_{yx}^\infty = 0  \mbox{ and } \bpam_{xx}^\infty  + \bpam_{xy}^\infty  < \bpam_{yy}^\infty \}, \mbox{ which implies}\\
    \cM^\complement &= \{ \mathbf{e}: \bpam_{yx}^\infty  = 0 \} \cap \{ \mathbf{e}:  \bpam_{yx}^\infty > 0  \mbox{ or } \bpam_{xx}^\infty  + \bpam_{xy}^\infty  \geq \bpam_{yy}^\infty \} = \{ \mathbf{e}:  \bpam_{yx}^\infty = 0, \bpam_{xx}^\infty  + \bpam_{xy}^\infty  \geq \bpam_{yy}^\infty \}. \nonumber
\end{align}
Observe that the first and second sub-classes in $\cM$ consider double and single-sided attack, respectively (at the limit); both classes consider acquisition.
An important question for a BP with attack is regarding the survival of the individual types and co-survival. Corollary \ref{corollary_BPA} of Theorem \ref{thrm1} given later provides answers to such questions. Prior to that, the next theorem derives the asymptotic analysis of \eqref{eqn_beta_ode_simple} and also shows that this analysis is sufficient for analysis of \eqref{ODE_BPA} (see proof in Appendix \ref{proof_thrmBPA}).

\begin{theorem}\label{thrmBPA}
Assume \ref{k1} and \ref{k2}. Then, \ref{a3} holds for \eqref{ODE_BPA}. Further, we have:
    \begin{enumerate}[label=(\roman*)]
        \item For ODE \eqref{eqn_beta_ode_simple}, no interior $\bc \in (0, 1)$ is  an attractor, $\bstar = 1$ is always an attractor, but  $\bstar = 0$ is  an attractor only if $\mathbf{e} \in \cM$.
        
        Further, again for \eqref{eqn_beta_ode_simple} in $[0,1]$: if $\mathbf{e} \in \cM$, then, $\bstar_r$, the unique zero of $g_\beta$, is the only repeller; while if $\mathbf{e} \notin \cM$, then $0$ is the only repeller.
        \item The attractors and repellers of ODE \eqref{eqn_beta_ode_simple} determine the attractor ($\cA$) and saddle ($\cR$) sets of ODE \eqref{ODE_BPA} respectively:
            \[
            \cA = 
            \begin{cases}
            \{ \mathbf{h}(1), \mathbf{h}(0)\}, &\mbox{ if } \mathbf{e} \in \cM,\\
            \{\mathbf{h}(1)\}, &\mbox{ if } \mathbf{e} \notin \cM,
            \end{cases}
            \mbox{ and }
            \cR = 
            \begin{cases}
            \{ \mathbf{0}, \mathbf{h}(\bstar_r)\}, &\mbox{ if } \mathbf{e} \in \cM,\\
            \{ \mathbf{0}, \mathbf{h}(0)\}, &\mbox{ if } \mathbf{e} \notin \cM, \ \ \mbox{ where}
            \end{cases}
            \]
            for example, $\mathbf{h}(1) = (\bpam^\infty_{xx}-1, \bpam_{xx}^\infty-1, \bpam_{xx}^\infty, \bpam_{xx}^\infty)$ and $\mathbf{h}(0) = (\bpam^\infty_{yy}-1, 0, \bpam^\infty_{yy}, 0)$.
        \item The combined domain of attraction of $\cA\cup \cR$, i.e., $\cD = \cD_I$ defined in \eqref{invariant_set}. \eop
    \end{enumerate}
\end{theorem}


\subsection{Analysis of random trajectory of BP with attack}\label{subsec_BPA}
By Theorem \ref{thrm1}, the following holds (proof in Appendix \ref{proof_cor_BPA}):
\begin{corollary}\label{corollary_BPA}
Consider the BP as in \eqref{eqn_dynamics_BPA}, and assume \ref{k1}-\ref{k2}. Then, we have:
\begin{enumerate}[label=(\roman*)]
    \item The assumption \ref{a3} holds for ODE \eqref{ODE_BPA}, and hence Theorem \ref{thrm1}(i) is applicable.
    \item The following is true w.p. $1$ for BP with attack:
    
    $\noindent$ $\bullet$ if $\mathbf{e} \in \cM$, either $\Ups_n$ converges to $\{\mathbf{0}, \mathbf{h}(0), \mathbf{h}(\bstar_r), \mathbf{h}(1)\}$ or  hovers around $\{\mathbf{0}, \mathbf{h}(\bstar_r)\}$, where $\bstar_r$ is as in Theorem \ref{thrm_attractors_beta} and
    
    $\noindent$ $\bullet$ if $\mathbf{e} \notin \cM$, either $\Ups_n$ converges to $\{\mathbf{0}, \mathbf{h}(0), \mathbf{h}(1)\}$ or hovers around $\{\mathbf{0}, \mathbf{h}(0)\}$. \eop
\end{enumerate}
\end{corollary}



Recall from Theorem \ref{thrmBPA}, ODE  for \eqref{ODE_BPA} has three types of saddle points: $\mathbf{h}(0)$ when $\mathbf{e} \notin \cM$, $\mathbf{h}(\bstar_r)$ when $\mathbf{e} \in \cM$ and vector $\mathbf{0}$ for all cases. 
The sample paths in which BP hovers around $\mathbf{0}$ or 
  $\mathbf{h}(0)$ or converges to/hovers around $\mathbf{h}(\bstar_r)$ indicate co-survival. Both populations survive in insignificant numbers in the first case, 
 $x$-population is comparatively small in the second case and both populations survive in large numbers in the last case. Further, only $x$ or $y$-population survives when the process converges to $\mathbf{h}(1)$ or $\mathbf{h}(0)$ respectively, see the end of Appendix \ref{appendix_prelim}.

We re-iterate that our approach does not provide the probability with which BP converges or hovers around different limit points of the ODE \eqref{ODE_BPA}.

\begin{figure}[htbp]
    \centering
    \includegraphics[trim = {0cm 1.5cm 0cm 0.4cm}, clip, scale = 0.35]{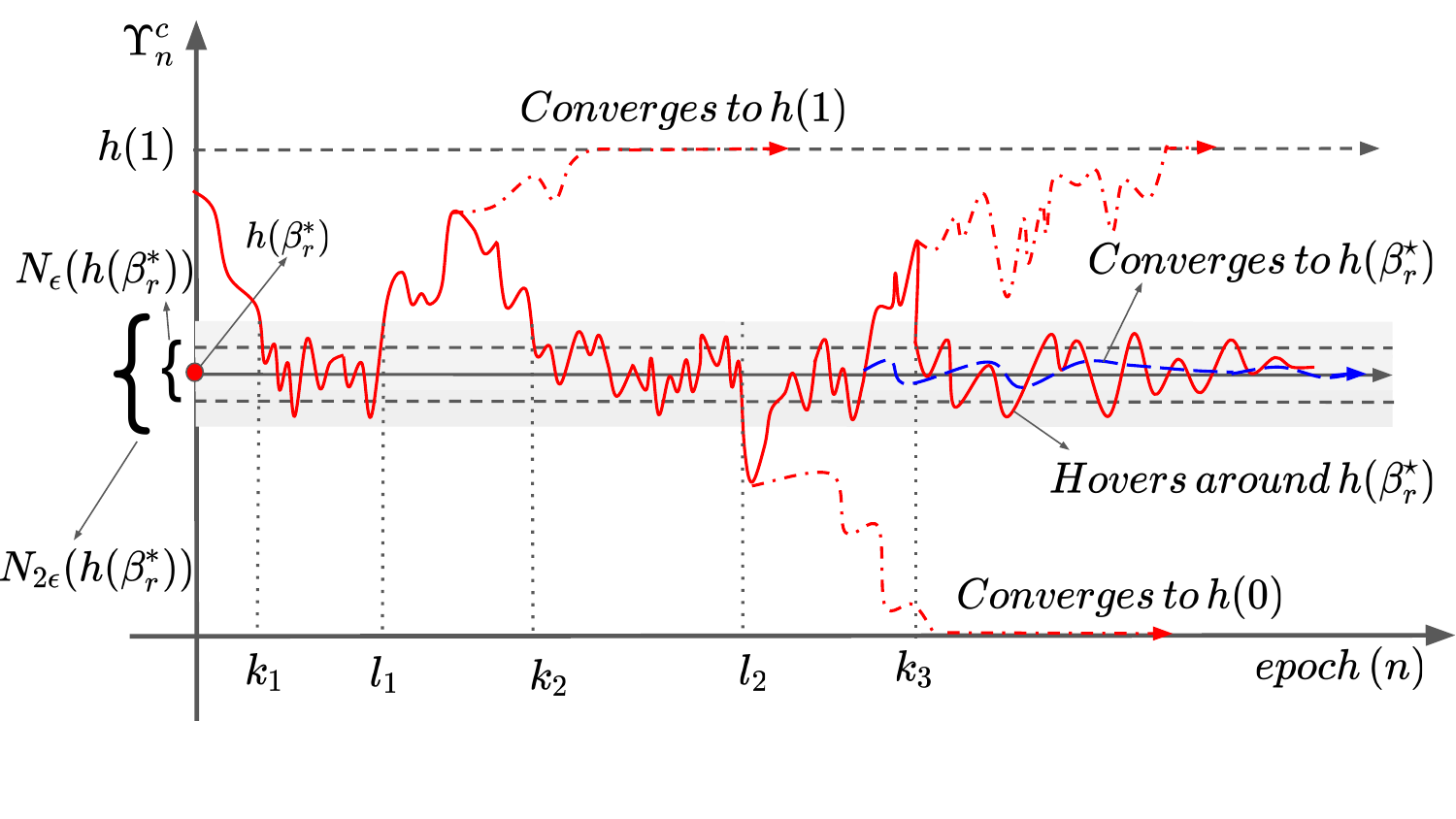}
    \vspace{-5mm}
    \caption{Behavior of BP with attack trajectory when $\mathbf{e} \in \cM$}
    \label{fig_repeller_BPA}
\end{figure}
Now, we would like to explain the behaviour of the BP with a pictorial representation in figure \ref{fig_repeller_BPA}. Consider $\mathbf{e} \in \cM$ and survival paths. Say, the process enters $\epsilon$-neighbourhood of $\mathbf{h}(\bstar_r)$ at epochs say $k_1, k_2, \dots$ (for some $\epsilon>0$), remains in its $2\epsilon$-neighbourhood for some epochs and then exits at epochs $l_1, l_2, \dots$ At every exit, it can either get attracted to $\mathbf{h}(0)$ or $\mathbf{h}(1)$ or it can re-enter the neighbourhood. The solid red line in the figure represents the sample path when the trajectory enters and exits the $\epsilon$-neighbourhood i.o., i.e., hovers around $\mathbf{h}(\bstar_r)$ with $\delta_1 = 2\epsilon$. Some sample paths can converge to $\mathbf{h}(\bstar_r)$ - see blue dashed line. Similar behaviour is exhibited when $\mathbf{e} \notin \cM$. Such hovering around is also observed in \cite{jagers2011population}, where switching between super-to-sub critical regimes occurs due to current population dependency.

\subsection{Application - Viral competing markets} \label{subsec_viral_competing}
In online social networks, content providers (CPs) share a variety of content, which is shared (again) by the recipients and thus may get viral (i.e., the number of copies of the post grows significantly with time). After reading the post, the user most likely loses interest in it forever. Thus, reading the post is analogous to death, while the number of new shares by a user is analogous to offspring. Further, unread and total (read $+$ unread) copies are analogous to the current and total population, respectively.

On such networks, contents often compete with each other (e.g., advertisements of similar products); when a new competing post (say $x$-type) is shared on the user's screen, the user might find $x$-post more attractive than an older $y$-post. This aspect leads to viral competing markets, where we say $x$-post has attacked and acquired the opportunities of $y$-post. Such attacks are dependent on the current copies. Further, the network is closed, and some users may share with previous recipients who would not be interested in the post again. Thus, the effective shares depend on the total copies. BP with attack precisely captures such dynamics (see \cite{agarwal2021co} for modelling details).

In \cite{agarwal2021co}, we analyzed such markets in a restricted setting, while Corollary \ref{corollary_BPA} can handle the generality mentioned here. Both the posts are prominent when the process converges to or hovers around $\mathbf{h}(\bstar_r)$. While, the convergence to $\mathbf{h}(0)$ or $\mathbf{h}(1)$ represents the dominance of one of the posts.

From Corollary \ref{corollary_BPA}, one can get more interesting insights. For instance, let $y$-CP be more influential, and thus $y$-post is shared more on average in the limit, so $\bpam_{xx}^\infty < \bpam_{yy}^\infty$. If the competition is ignored, the analysis is provided using independent BPs. Such analysis indicates the possibility of co-virality (both posts get viral simultaneously). However, when a typical user receives both posts, it may find $x$-post more appealing, leading to $\bpam_{xx}^\infty + \bpam_{xy}^\infty > \bpam_{yy}^\infty$ with $\bpam_{yx}^\infty = 0$. Therefore, $\mathbf{e} \notin \cM$, thus $\mathbf{h}(1)$ is a limit, which implies that $x$-post can dominate the post of more influential $y$-CP. Further, none of the limits indicate co-virality. 

On the other hand, when some users prefer the $y$-post ($\bpam_{yx}^\infty > 0$), while others prefer the $x$-post, then, co-virality is possible due to interior saddle point $\mathbf{h}(\bstar_r)$.

\section{Proportion-dependent BP and controlling fake news}\label{sec_prop_BP}

Consider a two-type (say $x$ and $y$ types) population-dependent BP. Here, each individual of any type lives for an exponentially distributed time, $\tau$. Before dying, it produces a random number of non-negative population-dependent offsprings of both $x$ and $y$ types. Now, the classical models assume the existence of the limit population-independent mean matrix, $M^\infty$ for the population-dependent mean matrix observed during transience. We consider the variant of BP, where in the limit, the mean matrix is proportion-dependent, and call such a BP as \textit{proportion-dependent BP} (PrD-BP). The dynamics of such a process can be captured by \eqref{evolve_x_up_time}.
    
    Observe that since $\lim_{t \to \infty} \bc(t)$ may not be known a priori, it is not easy to compute $M^\infty(\cdot)$ for PrD-BP. However, by virtue of Theorem \ref{thrm1} and upcoming Theorem \ref{thrmProp}, the possible limits of $M^\infty$ can be derived using the attractors of the ODE \eqref{eqn_beta_ode_simple} under 
 an additional assumption given below:
    
    
\begin{enumerate}[label=\textbf{A.\arabic*}, ref=\textbf{A.\arabic*}]
\setcounter{enumi}{4}
    \item  Assume $\minf_{ij}(\cdot)$, for $i, j \in \{x, y\}$ to be non-negative Lipschitz continuous on $[0, 1]$. Let $\cR_\beta$ be the finite set of repellers for the ODE \eqref{eqn_beta_ode_simple}. Further, let $\cA_\beta := \{\bstar_1, \dots, \bstar_k  \in [0, 1]\} $ be the attractor for the ODE \eqref{eqn_beta_ode_simple} for some $k < \infty$, with respective domain of attraction of each $\bstar_i$ as $\cD(\bstar_i)$ such that $\cup_{i=1}^k \cD(\bstar_i) = [0,1] - \cR_\beta$. \label{a5}
\end{enumerate}

\subsection{Analysis of ODE for PrD-BP}
We will now see that attractors/repellers of the ODE \eqref{eqn_beta_ode_simple} once again translates to appropriate attractors/saddle points for the ODE \eqref{eqn_ODE} as summarized below (see proof in Appendix \ref{appendix_B}):

\begin{theorem}\label{thrmProp}
Assume \ref{a5}. Then, the ODE \eqref{eqn_ODE} satisfies \ref{a3}. Further, the following statements are true for the ODE \eqref{eqn_ODE}:
\begin{enumerate}[label=(\roman*)]
    \item the attractor set, $\cA =  \{\mathbf{h}(\bstar) : \bstar \in \cA_\beta\}$, 
    \item the repeller set, $\cR = \{\mathbf{0}\}\cup \{\mathbf{h}(\bstar_r) : \bstar_r \in \cR_\beta\}$, and
    \item the combined DoA of $\cA$ and $\cR$, i.e., $\cD = \cD_I$. \eop
\end{enumerate} 
\end{theorem}

\subsection{Analysis of random trajectory of PrD-BP}
By Theorem \ref{thrm1}, we have: 

\begin{corollary}\label{corollary_prop}
Consider the dynamics as in \eqref{evolve_x_up_time}, and assume \ref{a1}-\ref{a2}, \ref{a5}. Then, following statements are true:
\begin{enumerate}[label=(\roman*)]
    \item the assumption \ref{a3} holds for the ODE \eqref{eqn_ODE}, and hence Theorem \ref{thrm1}(i) is true,
    \item $\Ups_n$ either converges to $\cA \cup \cR$, as $n \to \infty$ or hovers around $\cR$  w.p. $1$, where $\cA$ and $\cR$ are given in Theorem \ref{thrmProp}.  \eop
\end{enumerate}
\end{corollary}
Clearly, Corollary \ref{corollary_prop}(i) follows by standard results of ODE (e.g., \cite{piccini1984ordinary}) under \ref{a5}, and rest can be proved exactly as in Corollary \ref{corollary_BPA}.
We now discuss the implications of this Corollary.

$\bullet$ Let $M^\infty(\bc) \equiv M^\infty$, where all entries of constant matrix, $M^\infty$ are strictly positive, i.e., we have an irreducible mean matrix at limit which does not depend on population sizes/proportions. Then, by simple algebra (concavity/convexity of $\ga_\beta(\cdot)$), there exists a unique zero, $\bstar \in (0, 1)$, which is the only attractor for the ODE \eqref{eqn_beta_ode_simple} with DoA as $[0,1]$. 
Thus, $\Ups_n \to \{\mathbf{h}(\bstar), \mathbf{0}\}$, and therefore, the two populations co-survive in proportion equals to $\bstar$, if at all. This particular feature is analyzed already in the literature, except that we allow anything to occur during transience, along with total-population dependency (see Example $3$ in Section \ref{sec_6} where the process progresses as in BPA in the initial transition epochs, and then as PrD-BP).

$\bullet$ Now, consider $M^\infty(\bc) \equiv M^\infty$, where all entries of $M^\infty$ are strictly positive except say $\minf_{yx} = 0$, i.e., consider a decomposable proportion-independent process in the limit. Here, if $\minf_{xx} > \minf_{yy}$, the attractor for the ODE \eqref{eqn_beta_ode_simple} is uniquely given by:
\begin{align}\label{eqn_beta_decomposable}
    \bstar = \frac{\minf_{xx} - \minf_{yy}}{\minf_{xx} - \minf_{yy} + \minf_{xy} } \in (0,1).
\end{align}
 with DoA as $(0,1]$; thus, $1$ is not an attractor; further, $\bstar_r = 0$ is a repeller for the ODE \eqref{eqn_ODE}. In such a case, we have four kinds of sample paths: (i) both populations get extinct, (ii) only $x$-type gets extinct ($\Bc(t) \to 0$, in spite of $\bstar_r = 0$ being a repeller for the ODE), (iii) $x$-type population survives in small numbers (due to hovering around $\mathbf{h}(0)$), and (iv) the populations co-survive in proportion equals to $\bstar$. This result matches with \cite{jagers1969proportions} when one considers the Markovian variant of their age-dependent PnD model (see \eqref{eqn_beta_decomposable_pnd}).
 
 Otherwise if $\minf_{xx} \leq \minf_{yy}$, $\bstar = 0$ is the only attractor for the ODE with domain of attraction as $[0,1]$, which implies a.s. extinction of $x$-type population. This result matches with \cite[Theorem 1(ii)]{ranbir2019decomposable}, \cite{kesten1967limit} which considers population-independent decomposable BP (see Section \ref{desc_specialCases}).

$\bullet$ Lastly, if $M^\infty(\bc)$ is proportion dependent, then the set of limit proportions is $\cA_\beta \cup \cR_\beta$ (see \ref{a5}). This result is true even if the population dependent mean offspring functions, $m_{ij}(\om)$ depend only on $\bc$ (in the transience also). We discuss one example of the latter kind next directly while discussing the application.

\subsubsection{Controlling Fake News}\label{subsec_fake_news}
There is a variety of content propagating on OSNs, some of which are news articles. The news posts can be fake or authentic/real; fake news is fabricated (mis)information that propagates through OSNs like real news (see \cite{lazer2018science}). The consequences of fake news are fatal, and thus studies on the generation, propagation, detection, and control of fake news are welcome.

In \cite{kapsikar2020controlling} and \cite{kapsikar2020controllingarxiv}, we along with other authors focus on the control aspect of fake news. In particular, we design a control mechanism such that it can curb the propagation of fake news without significantly affecting the propagation of authentic news on OSNs; the users are asked to tag each post as authentic/fake. We model such controlled propagation as a BP; when a user receives the post, it tags the post as authentic ($y$-type) or fake ($x$-type). We further include the reluctance factor $\eta_c \leq 1$, which captures the hesitation/reluctance users have in forwarding the post after tagging it as fake. Such sharing results in more unread copies of either $x$-type and or $y$-type. The recipients tag the post based on the tag provided by the sender, their intrinsic ability to identify the actuality of the news, and the system-provided warning. 


In \cite{kapsikar2020controlling} and \cite{kapsikar2020controllingarxiv}, the system (OSN) exploits the 
users' collective wisdom by providing warning ($\omega(\bc)$) dependent  on the current proportion of fake tags ($\bc$). Thus, the offsprings/new shares depends on $\Bc(\cdot)$, according to the following mean matrix (see \cite{kapsikar2020controlling} and \cite{kapsikar2020controllingarxiv} for notations):
\begin{align*}
    M(\om) &= M(\bc) = M^\infty(\bc) = 
    \begin{bmatrix}
    \alpha_F^u\omega(\bc)m_\eta \eta_c & (1 - \alpha_F^u\omega(\bc))m_\eta \\
    \alpha_R^u\omega(\bc)m_\eta \eta_c & (1 - \alpha_R^u\omega(\bc))m_\eta
    \end{bmatrix}, \mbox{ where}\\
    \omega(\bc) &= \frac{w\bc}{\bc + b(1-\bc)} + \epsilon.
\end{align*}
From above mean matrix, it is clear that PrD-BP captures such propagation dynamics. In \cite{kapsikar2020controlling}, the analysis is available only for the case of  symmetric offsprings ($\eta_c = 1$). 
Using the results of this paper, one can extend the analysis for the case with $\eta_c < 1$; by  \cite[Lemma 3, Appendix]{kapsikar2020controllingarxiv}, the assumption \ref{a5} is proved, where we show that the unique attractor of the ODE \eqref{eqn_beta_ode_simple},  $\bstar \in (0,1)$, satisfies the following fixed point equation: 
\begin{align*}
    [\bstar \alpha^u_F \omega(\bstar) + (1-\bstar) \alpha^u_R \omega(\bstar)] [\eta_c + \bstar(1-\eta_c)] = \bstar.
\end{align*}
The domain of attraction for $\bstar$ is $[0,1]$, as $0 \mbox{ and } 1$ are repellers of the ODE \eqref{eqn_beta_ode_simple}, when $\epsilon > 0$.
One can analyze many other extensions of this process, for example,  with population dependent mean offsprings/shares, or with warnings which depend on $\ba$, as $\ba$ is indicative of fraction of users (so far) who tag the post as fake.

\section{Finite horizon approximation}\label{sec_6}
In Theorem \ref{thrm1}(i), we proved the finite time approximation of $\Ups_n$ using the autonomous ODE \eqref{eqn_ODE}; such an ODE is obtained using the limit proportion-dependent mean functions ($\minf_{ij}(\bc)$). However, directly using the population-dependent mean functions $m_{ij}(\om)$, one may anticipate better approximation in transience. 

We claim that ODE, $\dot{\ups} = \gna(\ups, t)$, constructed using the actual conditional expectation, $E[\mathbf{L}_n | \mathcal{F}_n] = \gna(\ups, t)$ given in \eqref{eqn_nonauto_g} better approximates the BP; recall, the difference term ${\cal E}_1^{n}(\cdot)$ 
 of \eqref{eqn_diff_term1} converges to $0$ as shown in the proof of Theorem \ref{thrm1}. The approximation should further improve when the new ODE is initialised with $\Ups_{n_m}$, and not with $\lim_{n_m \to \infty} \Ups_{n_m}$ as in Theorem \ref{thrm1}. From \eqref{eqn_nonauto_g}, the new ODE is non-autonomous and discontinuous. Also by \ref{a2},  the right hand side  $\gna(\Ups, t)$, converges to that of ODE \eqref{eqn_ODE}, $\ga(\Ups)$, as $t \to \infty$. Approximation by such non-autonomous ODE is proved for super-to-sub critical total population-dependent BP in \cite{agarwal2022saturated}.  

We support our claim using numerical examples.

\begin{example}
Consider a population-dependent BP with only one (say $x$-type) population, and let $\Cx(0) = 2$. Assume that initially, the population-dependent mean offsprings reduce linearly with an increase in total population size ($\ax$), and then gets fixed to $1.2$  as below:
\[
    m_{xx}(\om) =
    \begin{cases}
    3 - 0.002 \ax, \mbox{ if } \ax \leq 400,\\
    1.2, \mbox{ if } \ax > 400,
    \end{cases}
    \mbox{ for any } \om = (\cx, \ax).
\]
Clearly, the limit mean function is $\minf_{xx} = 1.2$, when $\cx \to \infty$. From FIGURE \ref{figure1}, one can see that the the curves ($\pc_n = \cx_n/n, \pa_n = \ax_n/n$ versus $n$, for all $n \geq n_m = 5$) for random trajectory (black curve) and non-autonomous ODE trajectory (red) are close by. However, the curve for autonomous ODE trajectory (blue) matches with the other curves only as $n$ grows large. It can also be seen from the plots that the random BP trajectory converges to the attractor of the autonomous ODE, as $n$ increases.
\begin{figure}[htbp]
\centering
\begin{minipage}{.5\textwidth}
  \centering
  \includegraphics[trim = {0cm 8cm 0cm 8cm}, clip, scale = 0.3]{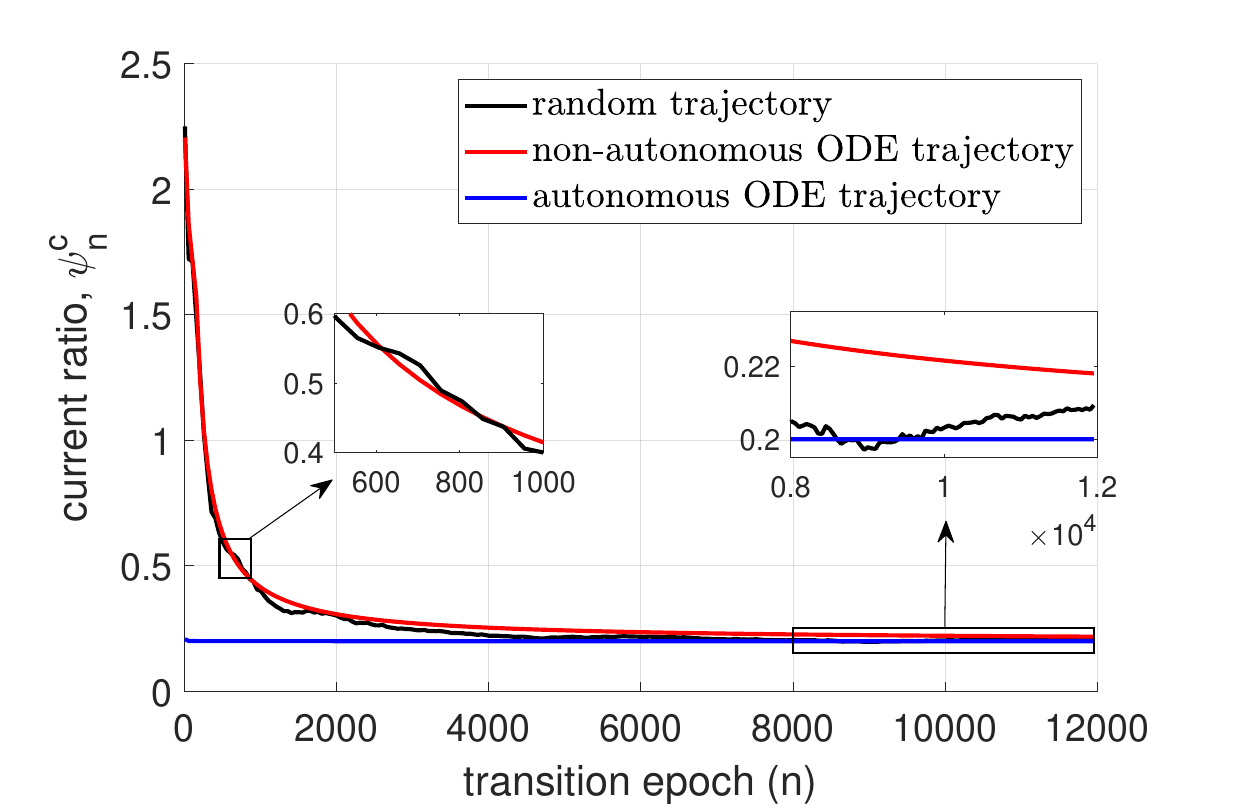} 
\end{minipage}%
\begin{minipage}{.5\textwidth}
  \centering
  \includegraphics[trim = {0cm 8cm 0cm 8cm}, clip, scale = 0.3]{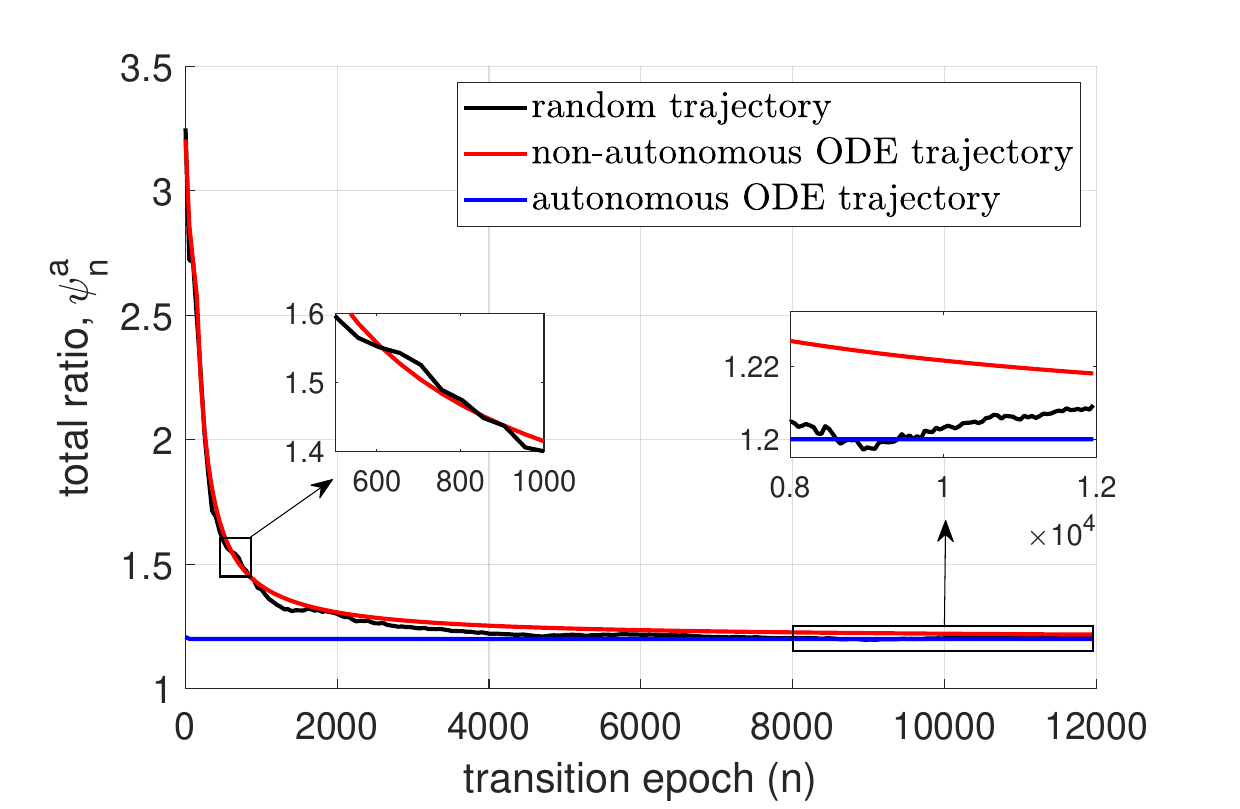}
\end{minipage}
\caption{Finite horizon approximation, single type PD-BP - one sample path}
\label{figure1}
\end{figure}
\end{example} 

\begin{example} Consider a PrD-BP with two population types ($x$ and $y$-type), $\Cx(0) = \Cy(0) = 100$ and the mean matrix:
\begin{align*}
M(\om) = M(\bc) = M^\infty(\bc) =
 \begin{bmatrix}
 6\bc & 2\bc\\
 4\bc & 5.6 \bc
 \end{bmatrix}.
 \end{align*}
\begin{figure}[htbp]
    \centering
    \includegraphics[trim = {0cm 6cm 0cm 6cm}, clip, scale = 0.4]{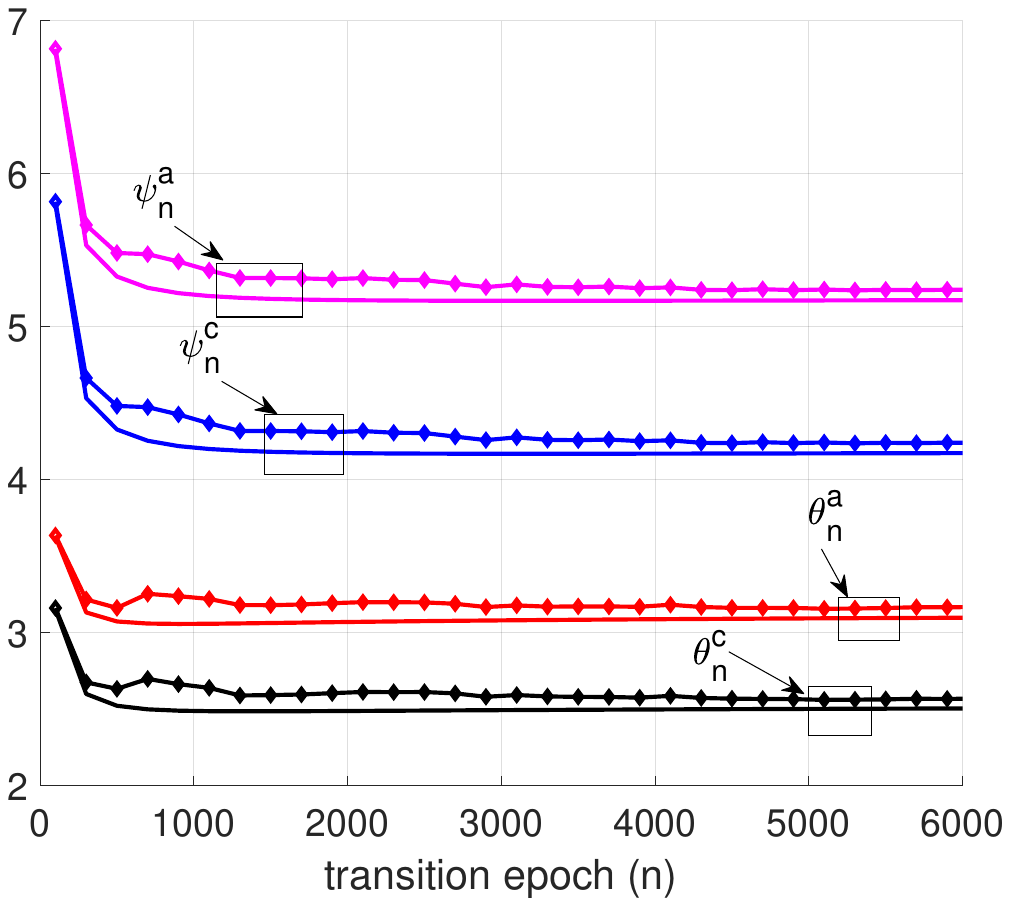}
    \caption{PrD-BP:  marked line- random trajectory (one sample path), and solid line- ODE trajectory}
\label{figure2}
\end{figure}
The above mean matrix is an instance of a variant of PrD-BP where the process is not in throughout super-critical regime, however, in a neighbourhood of the attractor of the corresponding ODE and stochastic system (in survival paths), the process is in super-critical regime. This paper does not cover the theoretical analysis of such processes, nonetheless, we numerically illustrate in FIGURE \ref{figure2} that the curves for random trajectory and ODE trajectory (for all $n \geq n_m = 100$) well match with each other; observe that the mean matrix has same structure from the start and hence the autonomous and non-autonomous ODE solutions are same, except for the initial values. We leave the analysis of such processes as a part of future work.
\end{example}

\begin{example} \label{eg3} Let $\Cx(0) = \Cy(0) = 1200$ and let the dynamics be as in BP with attack till $S^a$ is below a certain threshold, and then let the population progress with proportion-dependent mean offspring. Specifically, $M(\om) = M^t(\om)1_{\{s^a \leq 10^4\}} +  M^\infty(\bc)1_{\{s^a > 10^4\}}$, where
  \begin{align*}
M^t(\om) = 
\begin{bmatrix}
4 & -\min(2, \cy)\\
-\min(1, \cx) & 2.2
\end{bmatrix} \mbox{ and }  M^\infty(\bc) = \begin{bmatrix}
4 \bc + 1 & 9\bc+1\\
8\bc+1 & 2.2\bc+1
\end{bmatrix}.
\end{align*}
The process is in throughout super-critical regime. We plot one sample path of BP and corresponding solutions of autonomous and non-autonomous ODEs\footnote{The ODE trajectories are estimated using the well known Piccard's iterative method (e.g., \cite{piccini1984ordinary}).} (for all $n \geq n_m = 100$ and $T=12$). The current and total populations are in figure \ref{fig_eg3}, while the proportion $\bc(\ups_n)$ is in figure \ref{fig_eg3_beta}. From the plots, one can see that the non-autonomous ODE solution (dashed lines) better approximates the random BP trajectory (dotted lines), than the autonomous ODE (solid lines). As seen from the sub-figures, the non-autonomous ODE well captures the transition, unlike ODE \eqref{eqn_ODE}.

\begin{figure}[htbp]
\centering
\vspace{-4mm}
\begin{minipage}{.5\textwidth}
  \centering
  \includegraphics[trim = {0cm 8cm 0cm 8cm}, clip, scale = 0.32]{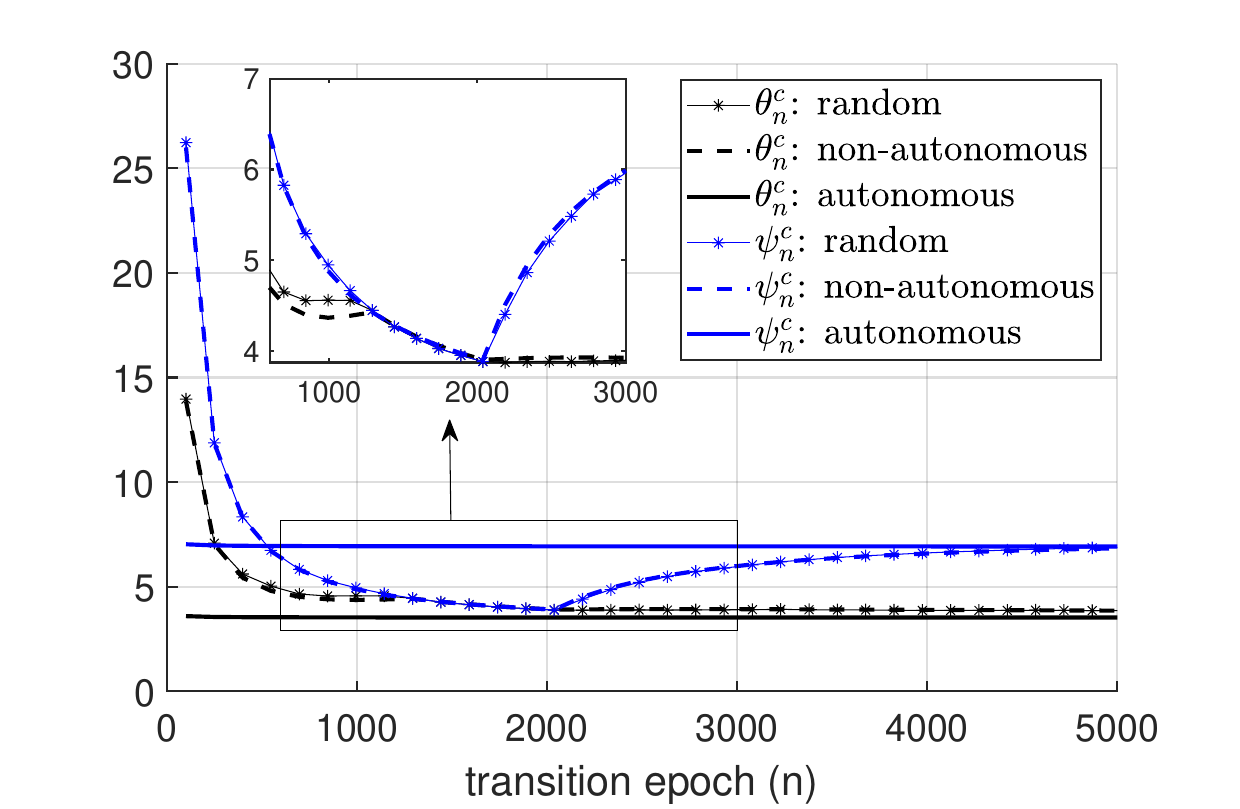} 
\end{minipage}%
\begin{minipage}{.5\textwidth}
  \centering
  \includegraphics[trim = {0cm 8cm 0cm 8cm}, clip, scale = 0.32]{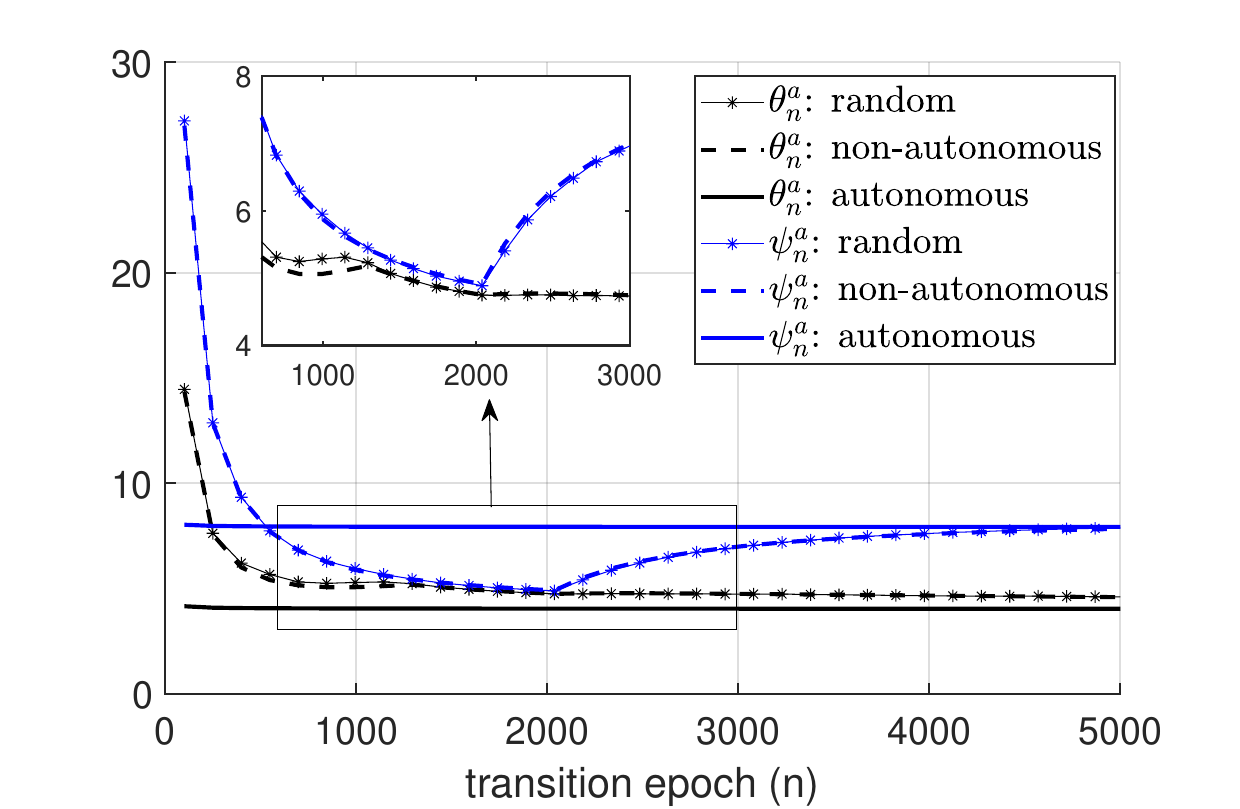}
\end{minipage}
\caption{Finite horizon approximation (current on left, and total on right side)}
\label{fig_eg3}
\end{figure}


\begin{figure}[htbp]
\centering
  \includegraphics[trim = {0cm 8.2cm 0cm 8.8cm}, clip, scale = 0.25]{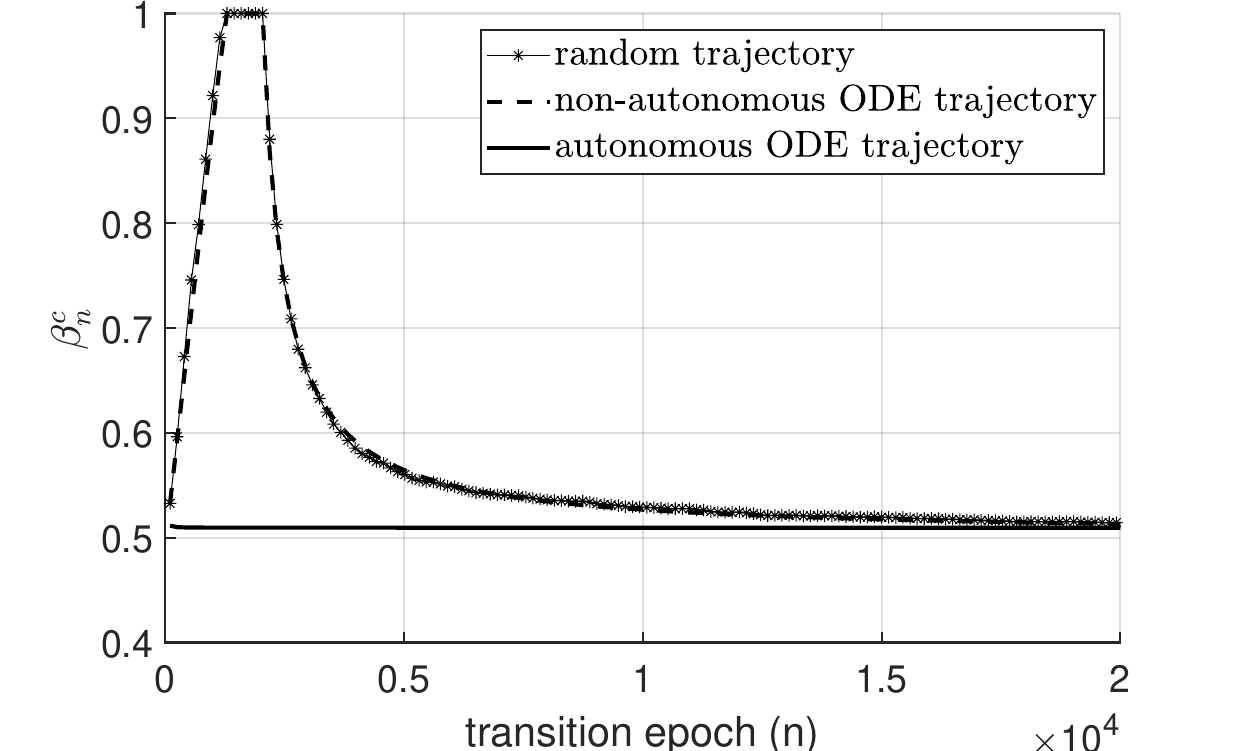}
  \caption{Proportion trajectory, $\bc_n$}
  \label{fig_eg3_beta}
\end{figure}
Initially, $x$-type individuals attack more aggressively than $y$-type and thus, the $y$-population depletes faster. In fact, by transition epoch ($1300$) proportion $\bc_n = 1$. Later, $M(\om) = M^\infty(\bc)$ does not have attack component, the $y$-population is regenerated and $\bc_n$ declines to $\approx 0.51$ indicating co-survival. This example also illustrates that the dynamics in transience (here, BP with attack) does not influence the limiting behaviour.

\end{example}

\subsection{Future extension: Super-to-sub critical BP}\label{subsec_supertosub}
The work in this paper so far has been about the BPs in the throughout super-critical regime. However, there is another important class of BPs, called super-to-sub critical BPs,  which we briefly discuss  here. 
\begin{definition}\label{defn_super_to_sub}
Let $d \geq 1$. A $d$-type population-dependent BP is said to be in \\ \underline{super-to-sub critical regime} if the maximal eigenvalue of mean matrix $M(\om) := [m_{ij}(\om)]_{ij}$ shifts from being $>1$ to $<1$ as $||\om^a||$ increases, where $\om^a := (a^i)_{1\leq i \leq d}$, and $||\cdot||$ is the Euclidean norm.   
\end{definition}
By the above definition, we want to consider scenarios where \textit{the population is initially in a super-critical regime and later (as time progresses) shifts to the sub-critical regime forever}. Such BPs are, for example, motivated by saturated content propagation in OSNs \cite{agarwal2022saturated}. We are able to capture such shift in regimes  due to the dependence on the total population ($||\Om^a||$), which is non-decreasing. 

In contrast, if one considers the dependence on non-monotonous current population, then the population fluctuates between super and sub-critical regimes. The authors in \cite{jagers2011population} capture such dynamics; they consider a continuous-time multi-type current population-dependent and age-dependent BP where the reproduction process shifts between super and sub-critical regimes, when the sum current population ($S^c(\cdot)$) goes below/above given carrying capacity ($K$). The authors show that either the population gets extinct, or the population size reaches the band around carrying capacity (i.e., $(K(1-\epsilon), K(1+\epsilon))$ for any $\epsilon > 0$) in a time of order $\log(K)$, and stays there for an exponentially long time ($\mathcal{O}(e^{cK})$, where $c > 0$ is some constant). 


In \cite{agarwal2022saturated}, we study \textit{saturated total-population dependent BP} (STP-BP), where a single population (say $x$-type) produces `total population-dependent offsprings' in Markovian and continuous-time framework. We show that $\Ax(t)$ converges and saturates to a limit as $t \to \infty$. Further, contrary to the known exponential growth in other existing BP models, with total population dependency as in Definition \ref{defn_super_to_sub}, $\Cx(t)$ grows exponentially initially and then declines to $0$.
We also prove a finite horizon result using SA-based methods; 
note that one can then derive the time analysis using the ODE trajectories, as in \cite{jagers2011population}. Furthermore, for a piece-wise linear total population-dependent mean function, we derive the (approximate) trajectories of current and total population sizes.

We believe that one can analyze a wide class of multi-type super-to-sub critical population-dependent-BPs (as in Definition \ref{defn_super_to_sub}) again using SA-based techniques, as described in this paper for throughout super-critical population-dependent-BPs, and as done in \cite{agarwal2022saturated} for STP-BP.

\section{Summary and conclusion}
We studied time-asymptotic proportion for a class of two-type continuous-time total-current population-dependent Markov BPs. We extended the stochastic approximation result to include the notion of hovering around the saddle points of an appropriate ODE and to analyze BPs. The summary to derive the limiting behaviour is:

(i) if the BP satisfies the assumption \ref{a1}, then the sum current population exhibits dichotomy with probability $1$ (see Lemma \ref{lemma_sum_pop});

(ii) identify the limit mean functions $\minf_{ij}(\bc)$ satisfying \ref{a2}, if required using the discussion in Appendix \ref{appendix_prelim} for BPs with negative  offspring or attack;

(iii) identify the attractors and repellers of one-dimensional ODE \eqref{eqn_beta_ode_simple};

(iv) identify the attractor and saddle sets of ODE \eqref{eqn_ODE} using (iii) and Theorem \ref{thrm_attractors_beta}; these provide the limit proportion; 

(v) Theorem \ref{thrm_attractors_beta} also facilitates the proof of \ref{a4} to conclude about limiting behaviour of BP via Theorem \ref{thrm1}.

Interestingly, the limit proportion of any BP depends only on the limit mean matrix, irrespective of the dynamics in transience. A finite-time approximation result is also provided. We analyzed a recently introduced variant of BP with attack and acquisition, and proportion-dependent BP under significantly more general conditions; such BPs respectively capture essential aspects of competing content propagation and fake news propagation over online social networks.

\begin{appendix}
\section{Some preliminary results}\label{appendix_prelim}
In this Appendix, we state some important auxiliary results, which are also helpful in further understanding of the subject at hand. The Lemma \ref{lemma_sum_pop} and the discussion thereafter provide insights into the derivation of the limit mean matrices of \ref{a2}.


\begin{lemma}{\bf [Dichotomy]}
\label{lemma_sum_pop}
Let $\underline{m} =: E[\underline{\offs}]$. Consider any BP satisfying \ref{a1}. Then, either $S^c(t) \to \infty$ exponentially at a rate at least $\lambda(\underline{m}-1)$, or $S^c(t) \to 0$ a.s.
\end{lemma}
\begin{proof}
Let $\Cx(0) = \cx_0$ and $\Cy(0) = \cy_0$. Consider a fictitious population-independent BP with single-type population, say $z$-type. Let $Z(0) = \cx_0 + \cy_0$. Each time an individual dies in the new process, assume that random number of offspring (distributed as $\underline{\offs}$ in \ref{a1}) are produced.  Further, assume that if $Z(t) = 0$ for some $t < \infty$, then, exactly $1$ individual is immigrated into the new system; this leads to the classical continuous time branching process with state-dependent immigration as in \cite{yamazato1975some}. Observe $\sum_{j=2}^\infty j P(\underline{\offs} = j) log(j) < \infty$ due to finite second moment assumption on $\overline{\offs}$ in \ref{a1}. Thus, by \cite[Theorems 6 and 8]{yamazato1975some}, $P(Z(t) \to \infty) = 1$, under \ref{a1}.

For completing the proof, we couple the embedded chains of the two BPs, for all $n \leq \nu_e$, where $\nu_e$ is the extinction epoch of the given system (see Section \ref{sec_probdesc_mainresult}); the offspring in the $Z(\cdot)$ branching process are given by $\underline{\offs}$ of \ref{a1}. If $\nu_e < \infty$, then $S_n^c = 0$ for all $n \geq \nu_e$. Otherwise, by coupling arguments, $S^c_n \geq Z_n$ for all $n$, and thus $S^c_n \to \infty$ as $n \to \infty$. Further, in the latter case, by \cite[Theorem 1, Chapter 1]{athreya2012classical}, the growth rate of $S_n^c$ is at least as large as that of $Z_n$, i.e., $\lambda(\underline{m}-1)$.
\end{proof}


\noindent \underline{\textbf{Limit mean matrices for BPs with negative offspring:}} 

In BPs with negative offspring, in the the survival sample-paths, by Lemma \ref{lemma_sum_pop}, $S^c_n \to \infty$. In such cases, one needs to identify the limit mean matrix of \ref{a2}. Say $0 < \liminf_{n \to \infty} \bc(\Ups_n) \leq \limsup_{n \to \infty} \bc(\Ups_n) <1$. Then, for such sample-paths, both populations would have exploded, i.e., $(\Cx_n, \Cy_n) \to (\infty, \infty)$. Hence, there are sufficient number of individuals to be attacked of both types, which results in the saturation of the number of attacks\footnote{To be realistic, the number of attacks by a single individual should saturate, i.e., for example, $\lim_{\cy \to \infty}m_{xy}(\cy) = \minf_{xy} < \infty$. The case with unsaturated attacks in easier to analyze, and one can easily prove for BP with attack that only one of the two population types survives with probability $1$.
\vspace{2.4mm}
}; thus, it is appropriate to consider $\minf_{xy}(\bc)$ as some constant for all $\bc \in (0,1)$, and so is the case with $\minf_{yx}(\bc)$.

On the other hand, say $\limsup_{n\to \infty} \bc(\Ups_n) = 1$, then, $\bc(\Ups_n) = 1$ i.o. This implies $\bc(\Ups_n) = 1$ for all $n$ large enough, as $\bc(\Ups_n) = 1$ is an absorbing state for processes with attack, like BP with attack and prey-predator BP. Thus, clearly $\minf_{xy}(\bc) = 0$ for $\bc = 1$. Similarly, $\minf_{yx}(\bc) = 0$ for $\bc = 0$.

\section{} \label{appendix_B}

Throughout the Appendix, we will consider the solution of the integral operator as the extended solution of ODE \eqref{eqn_ODE}. The fact that these two solutions are equivalent, is proved towards the end of the proof of Theorem \ref{thrm1}(i).

\noindent \textbf{Proof of Lemma \ref{lemma_equi_cont_thrm1} (contd.).} \label{proof_lemma1}
By \eqref{eqn_bounded_iterates}, $(\Ups^n(0))_n$ is bounded. We will now prove \eqref{eqn_footnote} for $(\Theta^{n, c}(t))$; it can be proved analogously for other components of $\Ups^n(\cdot)$. 
\hide{Towards this, re-write the scheme for $(\Theta_n^c)$ as follows (see
\eqref{eqn_bias},  \eqref{eqn_nonauto_ODE}):
\begin{align*}
    \Tc_{n+1} &= \Tc_n + \epsilon_n L_n^{\theta, c} =  \Tc_n + \epsilon_n (\delta M_n^{\theta, c} + g_\theta^{c}(\Ups_n) + D_n^{\theta, c}), \mbox{ where}
\end{align*}
$\delta M_n^{\theta, c} :=  L_n^{\theta, c} - g_\theta^{c}(\Ups_n) - D_n^{\theta, c} = L_n^{\theta, c} - \rho_\theta^{c}(\Ups_n, t_n)$. }
Observe from \eqref{eqn_diff_term1} and \eqref{eqn_diff_term2} that the interpolated trajectory can be re-written as:

\vspace{-4mm}
{\small
\begin{align}\label{eqn_interpolated_traj}
\begin{aligned}
    \Theta^{n, c}(t) &:= \Tc_n + \int_0^t g_\theta^c(\Ups^n(s)) ds +  \sum_{i=n}^{\eta(t_n+t)-1} \epsilon_i L_i^{\theta, c} - \int_0^t g_\theta^c(\Ups^n(s)) ds\\
    &=  \Tc_n + \int_0^t g_\theta^{c}(\Ups^n) ds + M^{n, \theta, c}(t) + \rho^{n, \theta, c}(t) + D^{n, \theta, c}(t), \mbox{ where}\\
    M^{n, \theta, c}(t) &:=  \sum_{i=n}^{\eta(t_n + t)-1}  \epsilon_i \left(L_i^{\theta, c} - \rho_\theta^{c}(\Ups_i, t_i)\right), \\ \rho^{n, \theta, c}(t) &:=  \sum_{i=n}^{\eta(t_n + t)-1}\epsilon_i g_\theta^{c}(\Ups_i) - \int_0^t g_\theta^{c}(\Ups^n) ds, \mbox{ and}  \\
    D^{n, \theta, c}(t) &:= \sum_{i=n}^{\eta(t_n + t)-1} \epsilon_i \left(\rho_\theta^{c}(\Ups_i, t_i) - g_\theta^{c}(\Ups_i)\right).
\end{aligned}
\end{align}}
\hide{Then, \eqref{eqn_interpolated_traj_1} can be re-written as: 

\vspace{-6mm}
{\small 
\begin{align}\label{eqn_interpolated_traj}
    \Theta^{n, c}(t) &:= \Tc_n + \sum_{i=n}^{\eta(t_n + t)-1}\epsilon_i (\delta M_i^{\theta, c} + g_\theta^{c}(\Ups_i) + D_i^{\theta, c}) \nonumber \\
    &= \Tc_n + \int_0^t g_\theta^{c}(\Ups^n) ds + M^{n, \theta, c}(t) + \rho^{n, \theta, c}(t) + D^{n, \theta, c}(t), \mbox{ where} \\
M^{n, \theta, c}(t) &:=  \sum_{i=n}^{\eta(t_n + t)-1}  \epsilon_i \delta M_i^{\theta, c}, \ \ \   D^{n, \theta, c}(t) := \sum_{i=n}^{\eta(t_n + t)-1} \epsilon_i D_i^{\theta, c}, \mbox{ and} \nonumber \\
\rho^{n, \theta, c}(t) &:=  \sum_{i=n}^{\eta(t_n + t)-1}\epsilon_i g_\theta^{c}(\Ups_i) - \int_0^t g_\theta^{c}(\Ups^n) ds. \nonumber 
\end{align}}}
Now, fix $T > 0$ and define the set $S^\delta_T := \{(s, t) : 0\leq t-s\leq \delta, 0\leq t \leq T \}$. Then:
\begin{align}\label{eqn_sup_bound} 
    \sup_{S_T^\delta} |\Theta^{n, c}(t) - \Theta^{n, c}(s)| 
     &\leq \sup_{S_T^\delta} \left|\int_s^t g_\theta^{c}(\Ups^n) dr\right| +  \sup_{S_T^\delta} \left|M^{n, \theta, c}(t) - M^{n, \theta, c}(s)\right| \nonumber \\
     &+ \sup_{S_T^\delta} \left|\rho^{n, \theta, c}(t) - \rho^{n, \theta, c}(s) \right| + \sup_{S_T^\delta} \left|D^{n, \theta, c}(t) - D^{n, \theta, c}(s) \right|.
\end{align}
To prove our claim, we begin with the first term  of   \eqref{eqn_sup_bound}. From  \eqref{eqn_ODE} and \eqref{eqn_bounded_iterates}, $|g_\theta^{c}(\Ups)| \leq \hat{m}$ for an appropriate $\hat{m} > 1$, for any $\Ups$, and, thus:

\vspace{-4mm}
{\small
\begin{align*}
     \left|\int_s^t g_\theta^{c}(\Ups^n) dr\right| 
     \hide{&\leq \int_s^t \hspace{-3mm} \left|g(\theta^{n, c}(r)) \right|dr \\
     &=  \int_s^t  \left|\bc\big(m_{xx}(\om) - 1\big) + (1-\bc) m_{yx}(\om) - \theta^{n, c}(r) \right|dr\\
     &\leq \int_s^t  \left|\bc (\hat{m} - 1) + (1-\bc) \hat{m}  \right|dr }
     \leq \hat{m}(t-s), \mbox{ so, }  \sup_{S_T^\delta} \int_s^t  \left|g_\theta^{c}(\Ups^n) \right|dr \leq  \delta \hat{m}.
\end{align*}}For the second term of   \eqref{eqn_sup_bound}, define $M_n^{\theta, c} := \sum_{i=0}^{n-1}\epsilon_i \left(L_i^{\theta, c} - \rho_\theta^{c}(\Ups_i, t_i)\right)$. Then, it is easy to prove that $(M_n^{\theta, c})$ is a Martingale   with respect to $(\mathcal{F}_n)$. Thus, using Martingale inequality, 
for each $\mu > 0$ (where, $E_n(\cdot)$ denotes the expectation conditioned on $(\mathcal{F}_n)$): 
$$
P\left\{\sup_{m\leq j \leq n} |M_j^{\theta, c} - M_m^{\theta, c}| \geq \mu \right\} \leq \frac{E_n\left|\sum_{i=m}^{n-1} \epsilon_i \left(L_i^{\theta, c} - \rho_\theta^{c}(\Ups_i, t_i)\right) \right|^2}{\mu^2}.
$$
Observe, $E\left[\left(L_i^{\theta, c} - \rho_\theta^{c}(\Ups_i, t_i)\right)\left(L_j^{\theta, c} - \rho_\theta^{c}(\Ups_j, t_j)\right)\right] 
= 0$ for $i < j$. Using this: 

\vspace{-2mm}
{\small
\begin{align*}
    P\left\{\sup_{m\leq j \leq n} |M_j^{\theta, c} - M_m^{\theta, c}| \geq \mu \right\} &\leq 
    \frac{\sum_{i=m}^{n-1} \epsilon_i^2 E_n\left| L_i^{\theta, c} - \rho_\theta^{c}(\Ups_i, t_i) \right|^2}{\mu^2}.
\end{align*}}
Note that under \ref{a1} and \eqref{eqn_bounded_iterates}, for some $K > 0$:
\begin{align*}
     \sup_n E_n|L_n^{\theta, c}- \rho_\theta^c(\Ups_i, t_i)|^2 &\leq
     \DetailK{\sup_n E\left|H_{n}\left(\offs_{xx, n}(\Om_{n-1}) - 1\right) + \overline{H}_{n} \offs_{yx, n}(\Om_{n-1}) - \Tc_{n-1}\right|^2\\
     &\leq} \sup_n E_n\left(\overline{\offs}_n - 1 \right)^2 + \sup_n E_n|\rho_\theta^c(\Ups_i, t_i)|^2 < K.
\end{align*}Thus, for every $n \geq m$:
\begin{align*}
    P\left\{\sup_{m\leq j \leq n} |M_j^{\theta, c} - M_m^{\theta, c}| \geq \mu \right\} &\leq  \frac{K}{\mu^2} \sum_{i=m}^{\infty} \epsilon_i^2.
\end{align*}
By first letting $n \to \infty$ (and using continuity of probability), then, letting $m \to \infty$, 
\begin{align}
    \lim_{m \to \infty} P\left\{\sup_{m\leq j } |M_j^{\theta, c} - M_m^{\theta, c}| \geq \mu \right\} &= 0 \mbox{ for each }\mu > 0.\label{eqn_equi_cont_M}
\end{align}
Now, by \eqref{eqn_equi_cont_M} and continuity of probability, for each $\mu > 0$:
\begin{align}\label{eqn_lim_sup_second_term}
 P\left\{\lim_{m \to \infty} \sup_{m\leq j } |M_j^{\theta, c} - M_m^{\theta, c}| \geq \mu \right\} = 0. 
\end{align}
Let $A_k := \lim_{m \to \infty} \sup_{m\leq j } |M_j^{\theta, c} - M_m^{\theta, c}| < 1/k$, then, $P(A_k) = 1$ for each $k > 0$. We further restrict our attention to  sample paths  $\omega \notin \underline{N :=  (\cap_k A_k)^c \cup \{\overline{\Pi} \nto \overline{m} \}}$. Now, the second term in \eqref{eqn_sup_bound} is upper bounded by $2\sup_{t\geq 0}|M^{n, \theta, c}(t)|$. For any $\omega \notin N$:
\hide{
\vspace{-4mm}
{\small
\begin{align*}
|M^{n, \theta, c}(t)| &= \left|\sum_{i=n}^{\eta(t_n + t)-1}\epsilon_i \delta M_i^{\theta, c} \right| 
= \bigg|M^{\theta, c}_{\eta(t_n + t)} - M^{\theta, c}_n\bigg|.
\end{align*}}
This gives us:}
\begin{align*}
\sup_{t\geq 0}|M^{n, \theta, c}(t)| 
&= \sup_{t \geq 0} |M^{\theta, c}_{\eta(t_n + t)} - M^{\theta, c}_n|  = \sup_{j \geq n} |M^{\theta, c}_{j} - M^{\theta, c}_n|
\end{align*}
\begin{align*}
\implies & \lim_{n \to \infty} \sup_{S_T^\delta}|M^{n, \theta, c}(t)|  \leq \lim_{n \to \infty} \sup_{\eta(t_n + t) \geq n} |M^{\theta, c}_{\eta(t_n + t)} - M^{\theta, c}_n| < 1/k,
\end{align*}where the last inequality holds because we have considered sample paths which are not in $N$.
Letting $k \to \infty$, we get, $M^{n, \theta, c}(\cdot) \to 0$ uniformly on each bounded interval. 

For the third term in \eqref{eqn_sup_bound}, observe that when  $t = t_k - t_n$ $(k > n)$, $\rho^{n, \theta, c}(t) = 0$. Thus, 
\hide{Towards this, for $k = n+1$, we have $t_k - t_n = \epsilon_n$:
\begin{align*}
    \rho^{n, \theta, c}(t_{n+1} - t_n) &=  \sum_{i=n}^{\eta(t_{n+1})-1}  \epsilon_i g(\tc_i) - \int_0^{\epsilon_n} g(\theta^{n, c}(s)) ds\\
    &= \epsilon_n g(\tc_n) - \epsilon_n g(\tc_n) = 0,
\end{align*}
where the second equality holds true because for $0 \leq s \leq \epsilon_n$, $\theta^{n, c}(s) = \tc_n$.
Let the claim be true for $k = n+l$, $l > 0$, then for $ k = n+l+1$, we have:
\begin{align*}
    \rho^{n, \theta, c}(t_{n+l+1} - t_n) &= \rho^{n, \theta, c}(t_{n+l} - t_n)  + \epsilon_{n+l}g(\theta_{n+m}^c) - \int_{t_{n+m}-t_n}^{t_{n+m+1}-t_n} g(\theta^{n, c}(s))ds = 0.
\end{align*}
Thus, by induction the claim is true. Next, we need to show that
$\rho^{n, \theta, c}(t) \to 0$ uniformly in $t$ as $n\to \infty$.} \hide{, i.e., for every $\epsilon > 0$, there exists $n_\epsilon$ such that for all $n \geq n_\epsilon$ and for all $t > 0$, $|\rho^{n, \theta, c}(t) - 0| < \epsilon$.}  for any $|t| \leq T$ (following similar steps as in first term, and noting $\epsilon_{\eta(t_n + t)} \le \epsilon_n$):
\begin{align*}
    |\rho^{n, \theta, c}(t)| &=
    \left| \int_{t_{\eta(t_n + t)} - t_n}^t g_\theta^c( \Ups^n) ds \right|
    \hide{\\
    &\leq  \int_{ t_{\eta(t_n + t)} - t_n}^t \left|g_\theta^{c, \infty}(\ups^n) ds \right| \leq (t - \eta(t_n + t)) \overline{m} }
    < \epsilon_n \hat{m}.
\end{align*}Thus, $\rho^{n, \theta, c}(\cdot)$ uniformly converges to $0$ as $n \to \infty$ on each bounded interval.  

For the last term in \eqref{eqn_sup_bound}, we claim that $D^{n,\theta, c}(t)$ also converges to $0$ uniformly on each bounded interval in $(0, \infty)$ as $n \to \infty$, for each $\omega \notin N$. Towards this, first consider $\omega \in N^c \cap \{S_n^c \to 0\}$, i.e, extinction paths. Then, $\rho_\theta^c(\Ups_i, t_i) = 0$ and $g_\theta^{c}(\Ups_i) = 0$ for all $i > \nu_e$. Thus, trivially $\lim_{n \to \infty}D^{n,\theta, c}(t) = 0$ for all $t \in (0, \infty)$.

Next, consider $\omega \in N^c \cap \{S_n^c \nto 0\}$; for such sample paths, we first derive a uniform positive lower bound for $\Pc_n$, required to prove the claim. To this end, analogous to $\overline{\Pi}_n$ defined in \eqref{eqn_overline_S_n}, one can define $\underline{\Pi}_n$ using $ \underline{\offs}$ given in \ref{a1}. Then, following similar steps as before, i.e., using strong law of large numbers and computing as in \eqref{eqn_bounded_iterates}, we get $\Pa_n \geq \Pc_n \geq \Delta$ for an appropriate $\Delta > 0$, for all $n\geq 1$. Thus, we have for each $i \geq 1$ (see $\tc$ component of \eqref{eqn_ODE}, \eqref{eqn_nonauto_g} and assumption \ref{a2}):

\vspace{-4mm}{\small
\begin{align*}
|D_i^{\theta, c}| = |\Bc_i(m_{xx}(\Om_i) -m_{xx}^\infty(\Bc_i)) + (1-\Bc_i)(m_{yx}(\Om_i) - m_{yx}^\infty(\Bc_i))| \leq \frac{2}{S_i^c} = \frac{2}{\Pc_i \eta(t_i)} \leq \frac{2}{\Delta i}. 
\end{align*}}This implies that, (recall $\epsilon_i = 1/(i+1)$)
\begin{align*} 
|D^{n, \theta, c}(t)| = \left|\sum_{i= n}^{\eta(t_n + t) - 1}\epsilon_i D_i^{\theta, c}\right| \leq \sum_{i= n}^{\eta(t_n + t) - 1} \frac{2}{\Delta i (i+1)} \leq \sum_{i= n}^{\infty} \frac{2}{\Delta i (i+1)}, \mbox{ for any }t.
\end{align*}Thus, $D^{n, \theta, c}(t)$ uniformly converges to $0$ as $n \to \infty$. In all, by \eqref{eqn_sup_bound} and above analysis, it is clear that for each $T>0$ and for any $\epsilon > 0$, there exists $n_\epsilon$ such that $\sup_{S_T^\delta} |\Theta^{n, c}(t) - \Theta^{n, c}(s)| < \epsilon $ for all $n \geq n_\epsilon$; hence $(\Theta^{n, c}(\cdot))$ is equicontinuous in extended sense. \eop

\noindent \textbf{Proof of Theorem \ref{thrm1} (ii).} \label{proof_thrm1}
The proof is constructed for sample paths $\omega \notin N$, however, for simplicity, we drop $\omega$ (see Lemma \ref{lemma_equi_cont_thrm1} for definition of set $N$). \hide{Consider the set $S(\omega)$ defined as:
$$
    S(\omega) := \{\ups: \Ups_n(\omega) \mbox{ exits } N_{\delta_{\mbox{\scriptsize{$\ups$}}}}(\ups) \mbox{ i.o. for some } \delta_\ups > 0\}.
$$
If $S(\omega)^c \cap \cR \neq \emptyset$
, then $\Ups_n(\omega) \to \cR$, and more precisely, $\Ups_n(\omega) \to \cR-S(\omega)$. Otherwise,  i.e., if $\cR \subseteq S(\omega)$, 
then $\Ups_n(\omega) \in \chi := \cD_b \cap \left(\cap_{\ups \in \cR} N_{\delta_\ups}^c(\ups)\right)$ i.o., observe $\chi$ is compact. For simpler notations, we drop $\omega$ henceforth. Therefore,}
By \ref{a4}, $\Ups_n\in \cD_b$ i.o. Since $\cD_b$ is compact, $({\Ups}_n)$ has a limit point ${\ups}_0 \in \cD_b$; then, there exists a sub-sequence $(n_k)$ such that ${\Ups}_{n_k} \to {\ups}_0$. Further, by (extended) equicontinuity of $(\Ups^n(\cdot))$, there exists further sub-sequence (denote it again by $(n_k)$, for simpler notations) $({\Ups}^{n_k}(\cdot))$ which converges to the extended solution $\ups(\cdot)$ of the ODE \eqref{eqn_ODE} uniformly on each bounded interval. Also observe, ${\Ups}^{n_k}(0) = {\Ups}_{n_k} \to {\ups}_0$, and recall $\ups(0) = \ups_0$  is the initial condition for ODE \eqref{eqn_ODE}. Under characterization of attractor or q-attractor in \ref{a4}, the ODE solution $\ups(t)$ converges to some $\ups^* \in (\cA \cup \cR)\cap \cD_I$ as $t \to \infty$. 

We will now show that for any $\delta_1 > 0$, $\Ups_n$ visits $N_{\delta_1}(\ups^*)$ i.o. We will also discuss other convergence aspects to complete the proof. Towards this, fix $\delta_1 > 0$.

{\bf Step A:} To begin with, assume $\ups^* \in \cA \cap \cD_I$.  Then, by \ref{a4} (local stability) it is possible to choose $0 < \delta_2 < \delta_1$ such that any ODE solution, ${\widetilde \ups}(\cdot)$, satisfies the following:
\begin{equation}\label{eqn_local_stability}
    {\widetilde \ups} (t) \in N_{\delta_1}(\ups^*)
    \mbox{ for all } t \geq 0, \mbox{ when initial condition } {\widetilde \ups}(0) \in cl( N_{\delta_2}(\ups^*) ).
\end{equation}
Further, by convergence of solution, $\ups(t) \to \ups^*$, thus there exists $T_{\delta_2} < \infty$ such that:
\begin{align}\label{eqn_dist_ODE_A}
    d({\ups}(t) , \ups^*) < \delta_2/2 \mbox{ for all } t \geq T_{\delta_2}.
\end{align}Now, following similar steps as in part (i) (see \eqref{eqn_dist_scheme_ODE_}), there exists $\overline{n} < \infty$ such that:
\hide{
Further, there exists $\overline{N}$ such that (by uniform convergence):
$$
\sup_{T_{\delta_2} \leq t \leq 2 T_{\delta_2}} d({\ups}^{n_k}(t) , {\ups}(t)) < \delta_2/2 \mbox{ for all } n_k \geq \overline{N}.
$$
Consider $t = t_l - t_{n_k}$ ($l > n_k$) such that $T_{\delta_2} \leq t \leq 2T_{\delta_2}$. Let $L:= \{l : T_{\delta_2} + t_{n_k} \leq t_l \leq 2T_{\delta_2} + t_{n_k}\}$. Then, we have (note that ${\ups}^{n_k}(t_l - t_{n_k}) = {\Theta}_l$):}
\begin{align}\label{eqn_dist_scheme_ODE}
\sup_{l \in L_k} d({\Ups}_l , {\ups}(t_l)) < \delta_2/2 \mbox{ for all } n_k \geq \overline{n},
\end{align}
for $L_k:= \{l : T_{\delta_2} + t_{n_k} \leq t_l \leq 2T_{\delta_2} + t_{n_k}\}$.
Using \eqref{eqn_dist_ODE_A} and \eqref{eqn_dist_scheme_ODE},  for all $n_k \geq \overline{n}$:
\begin{align}\label{eqn_dist_scheme_A}
\sup_{l \in L_k} d({\Ups}_l , \ups^*) &\leq \sup_{l \in L_k}d({\Ups}_l , {\ups}(t_l)) + \sup_{l \in L_k}d(\ups(t_l) , \ups^*) < \delta_2.
\end{align}Thus, ${\Ups}_n$ visits $N_{\delta_2}(\ups^*)$ i.o., and hence $N_{\delta_1}(\ups^*)$ i.o.

Henceforth, the proof is majorly as in proof of \cite[Theorem 2.3.1, pp. 39]{kushner2012stochastic}, except for \textit{few changes to consider convergence to q-attractors, not just attractors}. Contrary to the claim, assume that ${\Ups}_n$ exits $N_{\delta_1}(\ups^*)$ i.o. Thus, by \eqref{eqn_dist_scheme_A}, $\Ups_n$ moves from $N_{\delta_2}(\ups^*)$ to $\cD_b - N_{\delta_1}(\ups^*)$ i.o. 
Let $\overline{\Ups}^0(\cdot)$ be the usual linear interpolated trajectory of $\Ups_n$, i.e., 
\begin{align*}
    \overline{\Ups}^0(t_n) = \Ups_n, \mbox{ and } \overline{\Ups}^0(t) = \frac{t_{n+1}-t}{\epsilon_n}\Ups_n + \frac{t-t_{n}}{\epsilon_n}\Ups_{n+1} \mbox{ for } t \in (t_n, t_{n+1}).
\end{align*}
Then, there exists sequence $(l_j, r_j)$ such that (i) $\dots > r_j > l_j > r_{j-1} > l_{j-1}> \dots$, (ii) $r_j \to \infty$, (iii) $\overline{\Ups}^0(l_j) \in \partial N_{\delta_2}(\ups^*)$, $\overline{\Ups}^0(r_j) \in \partial N_{\delta_1}(\ups^*)$, and (iv) $\overline{\Ups}^0(t) \in cl(N_{\delta_1}(\ups^*)) -  N_{\delta_2}(\ups^*)$, for all $t\in [l_j, r_j]$. Consider the segments (one for each $j$) of $\overline{\Ups}^0(\cdot)$, i.e., consider  functions, $\q_j (t) := \overline{\Ups}^0(l_j+t)$ for any $ t\geq 0$;
observe by construction that for each $j,$ we have $\q_j (t) \in \left \{ \ups: \delta_2 <  d(\ups, \ups^*) \le \delta_1 \right \}$ for all $0 < t \le r_j - l_j$.


\textbf{Case (a):} Suppose there is a $T < \infty$ such that for some sub-sequence (call it $j$ again) $r_{j} - l_{j} \to T$. Now, consider a sub-sequence of $(\q_j(\cdot))$ which (again) converges to some solution of ODE, $\widetilde{\ups}(\cdot)$ uniformly over $[0,T]$.\footnote{The equicontinuity in extended sense can easily be extended to linear interpolated trajectories.} Then, $\widetilde{\ups}(0) \in \partial N_{\delta_2}(\ups^*)$ and $\widetilde{\ups}(T) \in \partial N_{\delta_1}(\ups^*)$. This contradicts \eqref{eqn_local_stability}. For $T = 0$, there is an obvious contradiction.

\textbf{Case (b):} If $r_{j} - l_{j} \to \infty$, then, $\widetilde{\ups}(0) \in \partial N_{\delta_2}(\ups^*)$ and $\widetilde{\ups}(t) \in cl(N_{\delta_1}(\ups^*)) -  N_{\delta_2}(\ups^*)$ for all $t > 0$. Then, it is a contradiction to $\ups^*$ being an attractor. 

In all, ${\Ups}_n \to \ups^*$; since $\ups^* \in \cA \cap \cD_I$ is arbitrary, we have ${\Ups}_n \to \cA \cap \cD_I$.

{\bf Step S:} 
Now consider $\ups^* \in \cR \cap \cD_I$. 
If $\nu_e < \infty$, i.e., in extinction sample paths, $\Ups_n \to {\bf 0}$ and we are done. For others, $\lim \inf_n \Pc_n > 0$ by Lemma \ref{lemma_sum_pop}.
Thus, with $\nu_e = \infty$ and   $\ups^* \in \cR \cap \cD_I$, 
 by Definition \ref{defn_q_as}, the initial condition $\ups_0  \in {\mathbb S} (\ups^*)$ with  $\bc(\ups_0) = \bc (\ups^*)$.

Similar to step A, by exponential stability (\ref{a4}), one can show that \eqref{eqn_local_stability} follows for any ODE solution $\widetilde{\ups}(\cdot)$ when initial condition  $\widetilde{\ups}(0) \in N_{\delta_2}(\ups^*) \cap \mathbb{S}(\ups^*)$. Further, clearly \eqref{eqn_dist_ODE_A}-\eqref{eqn_dist_scheme_A} also hold for this case. Thus, $\Ups_n$ visits $N_{\delta_1}(\ups^*) \cap \mathbb{S}(\ups^*)$ i.o. 

Further, if for every $\delta_1 > 0$, $\Ups_n$ does not exit $N_{\delta_1}(\ups^*) \cap \mathbb{S}(\ups^*)$ i.o., then $\Ups_n \to \Ups^* \in \cR\cap \cD_I$. Otherwise, for every $\delta_1 > 0$, $\Ups_n$ visits and exits $N_{\delta_1}(\ups^*) \cap \mathbb{S}(\ups^*)$ i.o. 
\eop




\hide{

\begin{lemma}\label{lemma_nonauto_auto_ode}
For any $T>0$,
\begin{align}\label{eqn_nonauto_auto_ode}
    g(\om(\widetilde{\ups}^n(s), t_n+s), \widetilde{\ups}^n(s)) \to g_\infty(\ups^*(s)) \mbox{ for all } s \in [0, T], \mbox{ as } n \to \infty.
\end{align}
\end{lemma}
\begin{proof} Consider the following, with $\widetilde{\ups}^n(s) = \om(\widetilde{\ups}^n(s), t_n+s)$:
\begin{align*}
    \left|g(\om(\widetilde{\ups}^n(s), t_n+s), \widetilde{\ups}^n(s)) - g_\infty(\ups^*(s))\right| &
    = \left|g(\widetilde{\ups}^n(s), t_n+s) - g_\infty(\ups^*(s))\right|\\
    &      \hspace{-5.5cm}  
    \leq \left|g(\widetilde{\ups}^n(s), t_n+s) - g(\ups^*(s), t_n+s)\right| + \left|g(\ups^*(s), t_n+s) - g_\infty(\ups^*(s))\right|
\end{align*}
Again by \textbf{A}.5, the second term converges to $0$, as $n \to \infty$. For the first term, we illustrate the claim for $\pa$, with $\widetilde{\om}_n(s) = \om(\widetilde{\ups}^n(s), t_n+s)$ and $\om_n^*(s) = \om(\ups^*(s), t_n+s)$:

\vspace{-2mm}
{\small
\begin{align*}
    \left|g_\psi^a(\widetilde{\ups}^n(s), t_n+s) - g_\psi^a(\ups^*(s), t_n+s)\right| 
    &\leq \bigg|\bc \sum_{j \in \{x, y\}} \bigg(m_{xj}(\widetilde{\om}_n(s)) - m_{xj}(\om^*_n(s))\bigg) \\
    &\hspace{-4.7cm}+ (1-\bc)\sum_{j \in \{x, y\}} \bigg(m_{yj}(\widetilde{\om}_n(s)) - m_{yj}(\om^*_n(s))\bigg)  - \left((\widetilde{\psi}^a)^n(s) - (\pa)^*(s)\right) \bigg|\\
    &\hspace{-4.7cm}\leq \sum_{i, j} |m_{ij}(\widetilde{\om}_n(s)) - \minf_{ij}(\bstar)| + \sum_{i, j} |m_{ij}(\om^*_n(s)) - \minf_{ij}(\bstar)| + |(\widetilde{\psi}^a)^n(s) - (\pa)^*(s)|\\
    &\hspace{-4.7cm}\leq \frac{4}{(\tilde{s}^{c, n})^\alpha} + \frac{4}{(s^{c, *})^\alpha} + |(\widetilde{\psi}^a)^n(s) - (\pa)^*(s)| \leq \frac{8}{(\Delta \eta(t_n+s))^\alpha} + |(\widetilde{\psi}^a)^n(s) - (\pa)^*(s)|.
\end{align*}}By uniform convergence of $(\widetilde{\ups}^n)$ on each bounded interval, and   $\eta(t_n + s) \to \infty$, as $n \to \infty$, the claim holds.
\end{proof}}

\vspace{2mm}
\newcommand{\OL}{\Omega}
\noindent \textbf{Proof of Theorem \ref{thrm_attractors_beta}.}\label{proof_thrm2}  Recall $\bc(\ups) := \tc/\pc$.  Consider the initial condition $\ups(0) \in \cD_I$ with $\pc(0) = 0$, then ODE \eqref{eqn_ODE} simplifies to $\dot{\ups} = -\ups$, which clearly has a unique solution and further $\ups(t) \to \mathbf{0}$ as $t \to \infty$. We claim that ${\mathbf{0}} \in \cR$ as we next show that   with $\pc(0) > 0$, the solution $\ups$ converges to other equilibrium points. 

Let $\pc(0) > 0$, and  say without loss of generality, $\bc(\ups(0)) \in {\cal N}_i^{-}$ for some $i$. By Lemma \ref{lemma_psi_c_general}, $\pc(t) > 0$ for all $t \geq 0$, thus ODE \eqref{eqn_ODE} simplifies to $\dot{\ups} = \mathbf{h}(\bc(\ups)) - \ups$. Consider the following smooth ODE, with initial condition $\ups(0) $ (by (c), the right hand side given below is Lipschitz continuous):
\begin{align}\label{eqn_func_fli}
\begin{aligned}
    \dot{\ups} &= f_l^i(\bc) - \ups, \mbox{ where } \\
    f_l^i(x) := \mathbf{h}(x)1_{\{x < x^*_i\}\cap N_i^*} &+ \mathbf{h}_l^* 1_{\{x \geq x^*_i\}} + \mathbf{h}_l^o 1_{x \leq \Delta_l^i}, \mbox{ with } \\
    \mathbf{h}_l^* := \lim_{x_n \up x^*_i} \mathbf{h}(x_n), \ \mathbf{h}_l^o := \lim_{x_n \downarrow \Delta_l^i} \mathbf{h}(x_n), &\mbox{ and } \Delta_l^i := \inf\{ \bc(\ups) : \bc(\ups) \in {\cal N}_i^*\}.
\end{aligned}
\end{align}Then, by \cite[Theorem 1, sub-section 1.4, pp. 6]{piccini1984ordinary}, the above smooth ODE has a unique solution, say $\ups^1(t)$. Let $\tau:= \inf\{ t : \bc(\ups^1(t)) = x^*_i\}$, then by Lemma \ref{lemma_tau_finite}, $\tau < \infty$. Observe that the solution of the original ODE \eqref{eqn_ODE}, with the same initial condition $\ups(0)$, coincides with $\ups^1(\cdot)$ for all $t < \tau$, as $\pc(t) > 0$ for all $t >0$ by Lemma \ref{lemma_psi_c_general} for such initial condition.  Now, let  $\ups^\tau := \ups^1 (\tau)$ and observe $\beta^c (\ups^\tau)  = x_i^*$. Using similar logic, one can prove that $x_i^*$ is an attractor for ODE \eqref{eqn_beta_ode_simple}. Further,  by uniqueness of the solutions of the smooth\footnote{The ODEs \eqref{eqn_func_fli} and \eqref{eqn_omega} are the two smooth ODEs.} ODEs, the solution of ODE \eqref{eqn_ODE} for $t > \tau$ is given by:
\begin{align}\label{eqn_ups2}
    \ups^2(t) = (\pc(t), x^*_i \pc(t), \pa(t), \ta(t)),
\end{align}where the three components of $\ups^2(\cdot)$, defined as $\OL(\cdot) := (\pc(\cdot), \pa(\cdot), \ta(\cdot))$ is the solution of the following initial value problem (IVP) for all $t \geq \tau$ (see \eqref{eqn_ODE}):
\begin{align}\label{eqn_omega}
\dot{\OL} = \mathbf{h}_i  - \OL, \mbox{ with } \OL (\tau) :=  \OL(\ups^*),  \mbox{ where constant, } \mathbf{h}_i := (h_{\psi}^c, h_{\psi}^a, h_{\theta}^a)|_{x^*_i}.
\end{align}Observe that $\bc(t) = x^*_i$ for all $t > \tau$ by (a). With this, $\ups(t) := \ups^1(t)1_{t < \tau} + \ups^2(t)1_{t > \tau}$ is the unique solution, which satisfies ODE \eqref{eqn_ODE}  for all $t \neq \tau$, and with initial condition $\ups(0)$. Thus, \eqref{eqn_ODE} satisfied \ref{a3}. Clearly from \eqref{eqn_omega}, 
$$
\ups(t) \to \mathbf{h}(x_i^*), \mbox{ where } \mathbf{h}(x_i^*) = (h^c_{\psi}, x_i^*h^c_{\psi}, h^a_{\psi}, h^a_{\theta})|_{x_i^*}.
$$
Similarly, one can show that $\ups(t) \to \mathbf{h}(x_i^*)$, if $\bc(\ups(0)) \in {\cal N}_i^{+}$. 

Thus, $\mathbf{h}(x_i^*)$ is an attractor for ODE \eqref{eqn_ODE}, with domain of attraction as $\cD_i := \{\ups\in \cD_I: \bc(\ups) \in  {\cal N}_i^*\} \cap \{\pc > 0\}.$ Since $x_i^* \in {\cal I}$ is arbitrary, $\cA = \{ \mathbf{h}(x_i^*) : x_i^* \in {\cal I}\}$, with corresponding domain of attraction as $\cD_{\cal A} =\cup_{1 \leq i \leq n} \cD_i$.  Also, ${\cal I}$ is an attractor for \eqref{eqn_beta_ode_simple}.

By hypothesis (b.i), any initial condition $\ups(0)$ with $\bc(\ups(0)) \in [0,1]-{\cal J}$ is already considered above. Now consider $\ups(0)$ with $\bc(\ups(0)) = y_i^* \in {\cal J}$, i.e., $\ups(0) \in {\mathbb S} (h(y_i^*))$. Then, the analysis follows as in \eqref{eqn_ups2}-\eqref{eqn_omega} to show that $\ups(t) \to \ups(y_i^*)$ as $t \to \infty$; the exponential convergence is clear from ODE \eqref{eqn_omega}. This proves that $\mathbf{h}(y_i^*)$ is a saddle point for ODE \eqref{eqn_ODE}.  Clearly, by (a), (b.ii)-(b.iii), $y_i^* \in {\cal J}$ is a saddle point for ODE \eqref{eqn_beta_ode_simple}. Hence, the theorem follows, as similar things are true for ${\mathbf 0}$.  \eop

\begin{lemma}\label{lemma_tau_finite}
The time $\tau$ defined in the proof of Theorem \ref{thrm_attractors_beta} is finite. 
\end{lemma}
\begin{proof}
By hypothesis (b),  $g_\beta(\cdot) > 0$ and continuous, for all $\bc \in {\cal N}_i^{-}$. Further,  $x_i^*$ is a point of discontinuity for $g_\beta$ and  $g_\beta(x_i^*) = 0$; thus $\bc(\mathbf{h}_l^*) = \lim_{x_n \up x_i^*} g_\beta (x_n) > 0$ (see \eqref{eqn_func_fli}), which implies, $\inf_{\{\bc \in  {\cal N}_i^{-}\}} g_\beta(\bc) > 0$. Observe $\tau$ is determined by $\bc$-component of $\ups^1(\cdot)$, the solution of  ODE \eqref{eqn_func_fli}. From \eqref{eqn_func_fli}, the latter is a continuous extension of the original ODE \eqref{eqn_ODE}, thus, the $\bc$-component of  the ODE \eqref{eqn_func_fli} can be uniformly lower bounded by $\inf_{\{\bc \in  {\cal N}_i^{-}\}} g_\beta(\bc) > 0$.
Thus, by Lemma \ref{lemma_general_ode_BPA}(a.ii), $\tau < \infty$.
\end{proof}

\begin{lemma}\label{lemma_general_ode_BPA}
Consider an initial value problem $\dot{z} = f(z, t)$, with $z(0) \in (z_0^l, z_0^u)$ where $f$ is a measurable function with finitely many discontinuities. 
\begin{enumerate}[label=(\alph*)]
    \item Say $f(z, t) > 0$, for all $z \in (z_0^l, z_0^u)$ and all $t$. Then:
    \begin{enumerate}[label=(\roman*)]
        \item  $z(\cdot)$ is an increasing function of $t$ till $\tau^u := \inf\{t : z(t) \geq z_0^u\}$.
        \item Say $f(z, t) > \delta$ for some $\delta > 0$, for all $z \in (z_0^l, z_0^u)$ and all $t$. Then, $\tau^u < \infty$.
    \end{enumerate}
    \item If $f(z, t) < 0$, for all $z \in (z_0^l, z_0^u)$ and all $t$, then $t \mapsto z(t)$ is a decreasing function till $\tau^l := \inf\{t : z(t) \leq z_0^l\}$, and if in addition $f(z, t) < -\delta$ for some $\delta > 0$, for all $z \in (z_0^l, z_0^u)$ and all $t$, then $\tau^l < \infty$.
\end{enumerate} 
\end{lemma}
\begin{proof}
We will provide the proof for part (a), and it can be done analogously for part (b). Contrary to the claim, let $\tau_1 < \tau_2 < \tau^u$ be two time points such that $z(\tau_1) \geq z(\tau_2)$, with  $z(\tau_1), z(\tau_2) \in (z_0^l, z_0^u)$. Then, we have:
$$
0 > z(\tau_2) - z(\tau_1) = \int_{\tau_1}^{\tau_2} f(z(s), s) ds,
$$
which is a contradiction to the hypothesis. 
Now if possible, let $\tau_u = \infty$, then $z(t) < z_0^u$ for all $t$ and $t \mapsto z(t)$ is an increasing function (as proved before).
Further, since $z(t) = z(0) + \int_{0}^t f(z(s), s) ds > z(0) + t \delta $, there exists $T_\delta > 0$ such that $z(t) \geq z_0^u$ for all $t \geq T_\delta$, which contradicts $\tau^u = \infty$. 
\end{proof}

\begin{lemma}\label{lemma_psi_c_general}
Let \ref{a2} and \ref{a3} hold. For any $0 < \epsilon < 2\leps  - 1$,  define $A_\epsilon := [2\leps  - 1-\epsilon, 2\ueps -1 + \epsilon]$. In case, $\pc(0) \in int(A_\epsilon)$ (interior) for some $\epsilon > 0$, then $\pc(t) \in A_\epsilon$ for all $t\geq 0$.  
Thus, if  $\pc(0) > 0$, then, $\pc(t) > \pc(0) - \delta$ for all $t \geq 0$ and for any $\delta > 0$.
\end{lemma}

\begin{proof} Recall from \eqref{eqn_ODE}, ODE for $\pc$ is $\dot{\pc} = h_\psi^c(\bc) 1_{\pc > 0} - \pc$. Now, one can lower bound $h_\psi^c(\bc) - \pc$ as, for all $t$ (by \ref{a1} and \eqref{eqn_bound_m_general}):
\begin{align}
    h_\psi^c(\bc) - \pc
    &\geq 2\bc \leps + 2\big(1 -\bc \big) \leps - 1 - \pc = 2\leps - 1- \pc.
\end{align}It is easy to observe (by Weierstrass Theorem) that there exists a strict positive uniform lower bound $l_I$ for any closed interval $I \subset (0, 2\leps-1)$ as below:
\begin{eqnarray}\label{eqn_g_psi_general}
\dot{\pc} \geq 2\leps-1- \pc \ge  l_I > 0 \mbox{ for any }  \pc \in I \mbox{ and for all } t. 
\end{eqnarray}For the first part, consider $I = [\pc(0), 2\leps-1-\frac{\epsilon}{2}]$, where $\pc(0) \notin A_\epsilon$. Then, by \eqref{eqn_g_psi_general}, $\dot{\pc} > l_I 
$ for all $\pc \in I$ and all $t$. By Lemma \ref{lemma_general_ode_BPA}(a), we have $\tau^u:=\inf\{t : \pc(t) \geq 2\leps-1-\frac{\epsilon}{2}\} < \infty$, i.e., $\pc(\cdot)$ enters $A_{\epsilon}$ from the left. 

We will now explicitly show that $\pc(\cdot)$ can not exit $A_\epsilon$, once it enters/starts in it (set $\tau^u = 0$ when $\pc(0) \in int(K_\epsilon)$). In contrast, say $\pc$ leaves $A_\epsilon$ and to the left.  Observe  $\pc(\tau^u) > 2\leps-1- \epsilon$. 
For $\pc$ to exit $A_\epsilon$, by continuity of $\pc$ (and Intermediate Value Theorem, IVT), there exist $2\leps-1-\epsilon < \underline{\nu}< \overline{\nu}< 2\leps-1$ such that for some $t_2 > t_1 > \tau^u$,  $\pc(t_2) = \underline{\nu}$ and $\pc(t_1) = \overline{\nu}$. Then, by MVT, we have:
$$
\dot{\pc}(s) = \frac{\pc(t_2) - \pc(t_1)}{t_2-t_1}  = \frac{\underline{\nu}-\overline{\nu}}{t_2-t_1} < 0, 
$$
for some $s \in (t_1, t_2)$. This is a contradiction as $\dot{\pc}(t) > 0$ for $\pc \in (0, 2\leps-1)$ and any $t$.
Conclusively, ODE solution $\pc(\cdot)$ enters $A_{\epsilon}$ from  left when $\pc(0) < 2\leps-1-\epsilon$, and does not exit $A_\epsilon$ from  left.

Similarly from \eqref{eqn_ODE}, $\mathbf{h}_\psi^c(\bc) - \pc$ 
can be  upper bounded  as (by \ref{a1} and \eqref{eqn_bound_m_general}):
\begin{align}
    \mathbf{h}_\psi^c(\bc) - \pc
    &\leq 2\bc \ueps + 2\big(1 -\bc\big)\ueps -1 - \pc = 2\ueps - 1 - \pc,
\end{align}and $\dot{\pc} \leq   2\ueps - 1 - \pc \le  u_I < 0$ for all $t$ and for any  $\pc \in I $ where $I \subset (2\ueps - 1, \infty)$ is any closed interval. Then, applying similar arguments as above, one can show that $\pc(\cdot)$ enters and does not exit $A_\epsilon$ from/to  right as well.
\end{proof}

\noindent \textbf{Proof of Theorem \ref{thrmBPA}.}\label{proof_thrmBPA} We first study ODE \eqref{eqn_beta_ode_simple}, using which we then analyze ODE \eqref{eqn_ODE}/\eqref{ODE_BPA}. Observe by definition of $\minf_{xy}(\cdot)$, $\minf_{yx}(\cdot)$ in \ref{k2} that $0, 1$ are equilibrium points of ODE \eqref{eqn_beta_ode_simple}. Further,  \textit{$g_\beta(\bc)$ is convex or concave in only $(0, 1)$, respectively if $\minf \leq 0$ or $\geq 0$}, as can be seen from below (see  \ref{k2} for definitions):
\begin{align}\label{eqn_g_beta}
\begin{aligned}
    g_\beta(\bc) &= \left(-\bpam_{yx}^\infty + \bc \widetilde{m}^\infty - (\bc)^2 \minf \right)1_{\bc \in (0, 1)}, \mbox{ where}\\
\widetilde{m}^\infty &:= \bpam^\infty _{xx} + \bpam_{xy}^\infty 
 -\bpam_{yy}^\infty + \bpam^\infty_{yx}, \mbox{ and } \minf := \bpam^\infty _{xx} - \bpam^\infty_{yy}.
\end{aligned}
\end{align}

At first by Lemma \ref{lemma_psi_c_general}, $\mathbf{0}$ is a saddle point for ODE \eqref{eqn_ODE} and hence for \eqref{ODE_BPA}.
Now, let $\minf \geq 0$, and consider the following two sub-cases.  

\textbf{Sub-case 1:} $\bpam_{xy}^\infty > 0$ and $\bpam_{yx}^\infty > 0$. Since $g_\beta(\cdot)$ is continuous in $(0, 1)$:
\begin{align}\label{eqn_case1}
    g_\beta(0^+) = \lim_{\delta \to 0} g_\beta(\delta) = -\bpam_{yx}^\infty  < 0, \mbox{ and }
    g_\beta(1^-) = \lim_{\delta \to 0} g_\beta(1-\delta)  = \bpam_{xy}^\infty  > 0.
\end{align}Therefore, there exists  a unique zero of $g_\beta$, say $\bstar_r \in (0, 1)$. Further by concavity,  $g_\beta(\bc) < 0$ when $\bc < \bstar_r$ and $g_\beta(\bc) > 0$ when $\bc > \bstar_r$. Thus, the result follows for this case by Theorem \ref{thrm_attractors_beta} with $x_1^* = 0$, $x_2^* = 1$ and $y^* = \beta^*_r$. That is, $\{0,1\}$ is the attractor set,  $\{\bstar_r\}$ is the repeller set for ODE \eqref{eqn_beta_ode_simple}. Thus, $\cA = \{\mathbf{h}(0), \mathbf{h}(1)\}$ is the attractor set and $\cD = \{\mathbf{0}, \mathbf{h}(\bstar_r)\}$ is the saddle set for ODE \eqref{eqn_ODE}, with combined domain of attraction, $\cD$ as in (v) of the Theorem.  

\textbf{Sub-case 2:} $\bpam_{xy}^\infty > 0$ and $\bpam_{yx}^\infty = 0$. Observe $\bpam^\infty_{xx} < \bpam^\infty_{yy}$ is not possible here, as it would contradict $\minf \geq 0$. Thus, $\bpam^\infty_{xx} \geq \bpam^\infty_{yy}$. Therefore, for any $\beta  \in (0, 1)$, $g_\beta(\beta) = \beta(1-\beta)(\bpam^\infty_{xx} - \bpam^\infty_{yy}) + \beta \bpam^\infty_{xy} > 0$. Further, $g_\beta(1^-) > 0$, as in case 1. Thus, the result follows for this case as well by Theorem \ref{thrm_attractors_beta} with $x_1^* = 1$ and $y^* = 0$. 



This completes parts (i) and (ii) for the case when $\minf \geq 0$. Analogously, one can prove (i) and (ii) when $\minf \leq 0$. Then, the proof is complete using Theorem \ref{thrm_attractors_beta}.   
\eop



\hide{
\noindent \textbf{Part (b):} Observe that $g_\beta(\cdot)$ with constant $m_{ij}^\infty(\cdot)$ can be re-written (as in \eqref{eqn_g_beta}):
\begin{align}
    g_\beta(\bc) = \minf_{yx} + \bc \widetilde{m}^\infty - (\bc)^2 \minf, \mbox{ for all } \bc \in [0,1].
\end{align}
Unlike part (a), $0, 1$ are not equilibrium points for the ODE \eqref{eqn_beta_ODE}. Now, as in part (a), we provide the proof for $\minf \geq 0$, which implies $g_\beta(\bc)$ is a concave function for all $\bc \in [0,1]$. 

\textbf{(i)} Consider $\minf_{xy} > 0$. At first, let $\minf_{yx} > 0$. Observe that $g_\beta(0) = \minf_{yx} > 0$ and $g_\beta(1) = -\minf_{xy} < 0$. Thus, there exists a unique $\bstar \in (0, 1)$ such that $g_\beta(\bstar) = 0$, $g_\beta(\bc) > 0$ for all $\bc \in [0, \bstar)$ and $g_\beta(\bc) < 0$ for all $\bc \in (\bstar, 1]$. 
This implies that the function $t \mapsto \bc(t)$ is strictly increasing (decreasing) till $\tau^u := \inf\{t : \bc(t) \geq \bstar\}$, if $\bc(t_0) \in [0, \bstar)$ for some $t_0 \in [0, \tau_u)$ (respectively, till $\tau^l := \inf\{t : \bc(t) \leq \bstar\}$, if $\bc(t_0) \in (\bstar, 1]$ for some $t_0 \in [0, \tau_l)$). As a result, $\bstar$ is AS with domain of attraction as $[0,1]$. 

Secondly, let $\minf_{yx} = 0$ with $\minf_{xx} > \minf_{yy}$. Here, $g_\beta(0) = 0$, however $g_\beta(0^+) > 0$. Thus, $\bstar = \widetilde{m}^\infty/\minf \in (0, 1)$ is AS with domain of attraction as $(0,1]$. 

\textbf{(ii)} Consider $\minf_{xy} > 0$ and $\minf_{yx} = 0$ with $\minf_{xx} \leq \minf_{yy}$. Note that $g_\beta(\bc) < 0$ for all $\bc \in (0, 1]$ and $g_\beta(0) = 0$. Then, as in (b.i), $\bstar = 0$ is AS with domain of attraction as $[0,1]$.}

\noindent \textbf{Proof of Corollary \ref{corollary_BPA}.}\label{proof_cor_BPA}
Given limit mean functions as in \ref{k2}, the assumption \ref{a3} is guaranteed by Theorem \ref{thrmBPA}. We now prove the assumption \ref{a4}.

$\cA$ and $\cR$ are the attractor and saddle sets of ODE \eqref{ODE_BPA} respectively, with subset of the combined domain of attraction as $\cD_I$, as identified in Theorem \ref{thrmBPA}. Towards getting a compact sub-domain of $\cD_I$, as in \eqref{eqn_pi_n}, 
 from \eqref{Eqn_XnYnSnetc_general}, \eqref{eqn_dynamics_BPA} and \ref{k1}, one can bound $\Pa_n$:
\begin{align*}
    \Pa_n \leq \overline{\Psi}_n^a :=  \frac{1}{n}\left(\sum_{k=1}^{\min\{\nu_e, n\}} \left( \overline{\xi}_{xx, k} + \overline{\xi}_{yy, k} \right)1_{\{\Pc_k > 0\}} + s_0^c \right).
\end{align*}
As before, $\overline{\Psi}_n^a \to E[\overline{\xi}_{xx,1} + \overline{\xi}_{yy,1}]$ a.s. in survival paths and $\overline{\Psi}_n^a \to 0$ in extinction paths, as $n \to \infty$. Thus, $\cS :=  \cD_I \cap  \left\{\ups :  \pa \in [0, E[\overline{\xi}_{xx,1} + \overline{\xi}_{yy,1}] ]\right\}$ is the compact subset of $\cD_I$ and $p_{b} := P(\Ups_n \mbox{ visits } \cS \mbox{ i.o.}) = 1$. Hence, by Theorem \ref{thrmBPA} and Theorem \ref{thrm1}(ii), we have $\Ups_n \to \cA \cup \cR$ with probability $1$. \eop

\end{appendix}

\bibliographystyle{imsart-number.bst} 

\hide{
\newpage
Below figure explains how to extend theorem 2 for $x_i^*$ which are continuity points.
\begin{figure}[htbp]
    \centering
    \includegraphics[scale = 0.07]{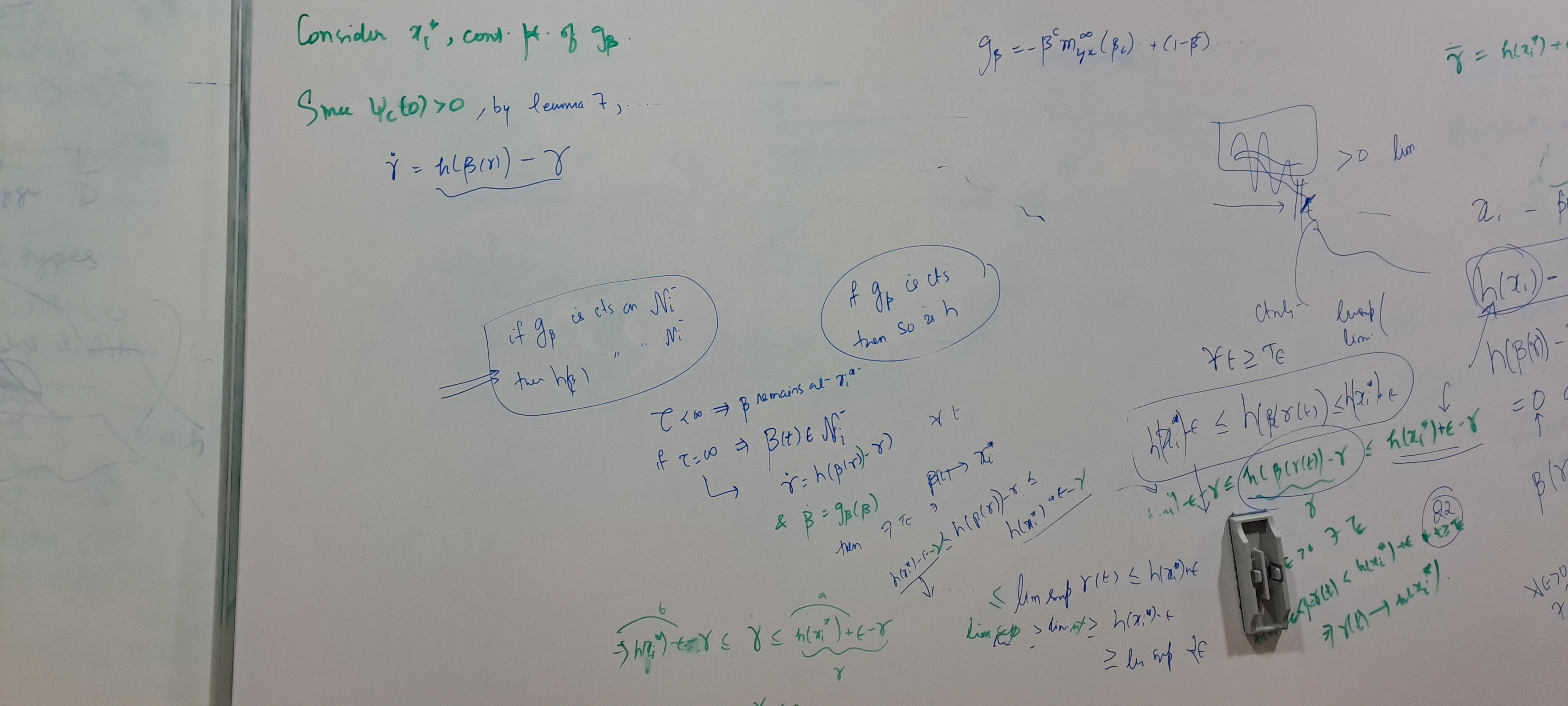}
\end{figure}

Below figure explains why definition of repellers is modified. Basically one can see that deviations in $\pc$-axis leads to convergence towards attractors. Else, $\ups$ converges to $0$ vector. 
\begin{figure}[htbp]
    \centering
    \includegraphics[scale = 0.07]{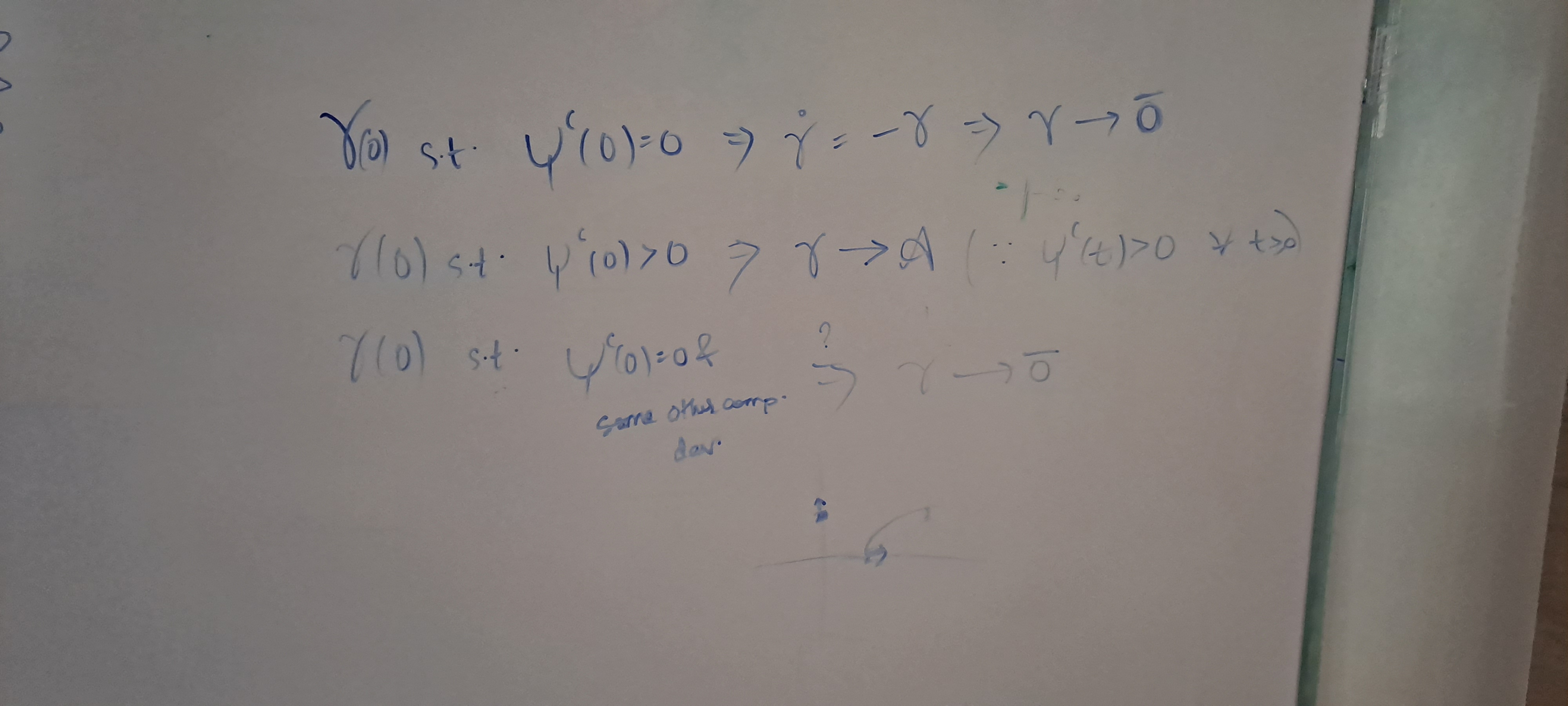}
\end{figure}

Below figure illustrates how $\cD_I$ is invariant set (see the commented line along definition).
\begin{figure}[htbp]
    \centering
    \includegraphics[scale = 0.07]{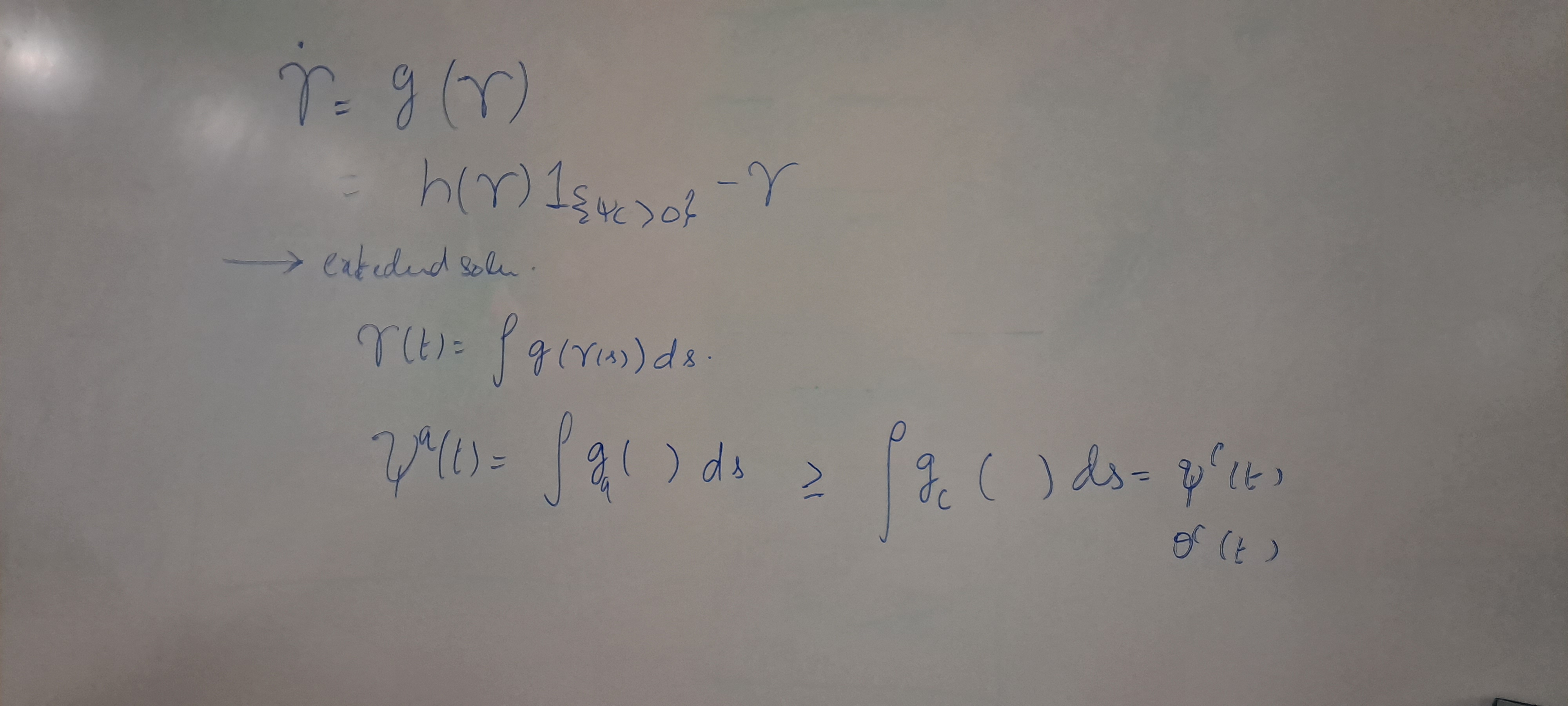}
\end{figure}}

\end{document}